\documentclass[12pt,letterpaper]{article}

\usepackage{multirow}%

\usepackage{dsfont}
\usepackage{mathrsfs}
\usepackage{amssymb,amsmath,amsfonts,amsthm}
\usepackage[parfill]{parskip}
\usepackage{enumitem} % for custom enumeration
\setlist[enumerate,1]{label={(\roman*)}}
\usepackage{url} % for URL link
\usepackage{xspace}
\usepackage[most]{tcolorbox}
\usepackage{scrextend}

\usepackage[ruled,vlined]{algorithm2e}

\usepackage{changepage}
\usepackage{thmtools}

\usepackage{graphicx}
\usepackage{subcaption}

\usepackage{tikz}
\usetikzlibrary{calc}
\usetikzlibrary{positioning, backgrounds, fit, patterns, patterns.meta}
\usetikzlibrary{decorations.pathmorphing}
\usetikzlibrary{intersections}

\usepackage{xcolor} % For custom colors

\definecolor{FigRed}{RGB}{046,037,133}
\definecolor{FigBlue}{RGB}{105,166,202}
\definecolor{FigGreen}{RGB}{148,203,236}
\definecolor{FigOrange}{RGB}{219,133,052}
\colorlet{FigPoolFill}{FigBlue!24}
\colorlet{FigCoreFill}{FigOrange!24}

\DeclareRobustCommand{\setcaptionbox}[2]{%
	\tcbox[on line, colback=#1, colframe=#1, boxrule=0pt,
		arc=2pt, boxsep=2pt, left=2pt, right=2pt,
		top=.3pt, bottom=.3pt]{\ensuremath{#2}}%
}

\tikzset{
	containment/ambient/.style={
		draw=black!65,
		fill=black!2,
		line width=.1pt,
		rounded corners=4pt,
		inner sep=0pt
	},
	containment/selected/.style={
		draw=FigRed,
		fill=FigGreen!22,
		line width=.1pt,
		rounded corners=3pt,
		inner sep=0pt
	},
	containment/core/.style={
		draw=FigBlue,
		fill=FigBlue!12,
		line width=.1pt,
		rounded corners=2pt,
		inner sep=0pt
	},
	containment/overlap/.style={
		pattern={Lines[angle=45, distance=24pt, line width=12pt]},
		pattern color=FigCoreFill,
		rounded corners=2pt,
		draw=none
	},
	containment/label/.style={
		font=\small,
		anchor=north west,
		inner sep=0pt
	}
}

\definecolor{linkGreen}{RGB}{175,0,0}

\definecolor{linkblue}{RGB}{0, 102, 204} 
\usepackage[
colorlinks=true,
linkcolor=linkblue,     % Color for links (e.g., theorems, equations)
citecolor=linkGreen,     % Color for citations
urlcolor=linkblue       % Color for URLs
]{hyperref}

\usepackage{orcidlink}

\usepackage{dirtytalk} % quotations with \say
\usepackage{mathtools} % to use \coloneqq

\usepackage[capitalize,nameinlink]{cleveref}
\DeclareMathSizes{7.3}{7.3}{5}{5}

\Crefname{figure}{Figure}{Figures}
\Crefname{claim}{Claim}{Claims}
\Crefname{algorithm}{Algorithm}{Algorithms}

\crefformat{equation}{\textup{#2(#1)#3}}
\crefrangeformat{equation}{\textup{#3(#1)#4--#5(#2)#6}}
\crefmultiformat{equation}{\textup{#2(#1)#3}}{ and \textup{#2(#1)#3}}
{, \textup{#2(#1)#3}}{, and \textup{#2(#1)#3}}
\crefrangemultiformat{equation}{\textup{#3(#1)#4--#5(#2)#6}}%
{ and \textup{#3(#1)#4--#5(#2)#6}}{, \textup{#3(#1)#4--#5(#2)#6}}{, and \textup{#3(#1)#4--#5(#2)#6}}

\Crefformat{equation}{Equation~\textup{#2(#1)#3}}
\Crefrangeformat{equation}{Equations~\textup{#3(#1)#4--#5(#2)#6}}
\Crefmultiformat{equation}{Equations~\textup{#2(#1)#3}}{ and \textup{#2(#1)#3}}
{, \textup{#2(#1)#3}}{, and \textup{#2(#1)#3}}
\Crefrangemultiformat{equation}{Equations~\textup{#3(#1)#4--#5(#2)#6}}%
{ and \textup{#3(#1)#4--#5(#2)#6}}{, \textup{#3(#1)#4--#5(#2)#6}}{, and \textup{#3(#1)#4--#5(#2)#6}}

\newtheorem{theorem}{Theorem}

\newtheorem{corollary}[theorem]{Corollary}
\newtheorem{observation}[theorem]{Observation}

\newtheorem{problem}[theorem]{Problem}
\newtheorem{lemma}[theorem]{Lemma}

\theoremstyle{definition}

\theoremstyle{remark}
\newtheorem{remark}[theorem]{Remark}
\newtheorem{claim}[theorem]{Claim}

\definecolor{mainresultcolor}{HTML}{F9E9D9}
\newenvironment{boxmainresult}{\begin{theorem}}{\end{theorem}}
\tcolorboxenvironment{boxmainresult}{
	colback=mainresultcolor,
	colframe=white,
	colbacktitle=mainresultcolor,
	coltitle=mainresultcolor
}

\newenvironment{pocd}[1]{\begin{proof}[Proof of {C}laim #1.]}{\end{proof}}

\renewcommand{\leq}{\leqslant}
\renewcommand{\geq}{\geqslant}
\renewcommand{\le}{\leqslant}
\renewcommand{\ge}{\geqslant}

\newcommand{\cC}{\ensuremath{\mathcal{C}}}
\newcommand{\cD}{\ensuremath{\mathcal{D}}}

\newcommand{\cF}{\ensuremath{\mathcal{F}}}

\newcommand{\cN}{\ensuremath{\mathcal{N}}}

\newcommand{\cP}{\ensuremath{\mathcal{P}}}

\newcommand{\bN}{\ensuremath{\mathbb{N}}}

\newcommand{\her}{\mathrm{Her}}

\newcommand{\rep}{\operatorname{rep}}
\newcommand{\defect}{\mathrm{defect}}
\newcommand{\good}{\mathrm{good}}
\newcommand{\pro}{\mathrm{pro}}
\newcommand{\tree}{\mathrm{tree}}
\newcommand{\ddef}{\mathrm{def}}
\newcommand{\prodef}{\mathrm{prodef}}

\usepackage{geometry}
\usepackage[utf8]{inputenc}
\usepackage{aliascnt}

\usepackage{amssymb,amsmath,amsthm,amsfonts,xcolor}

\makeatletter

\def\@gifnextchar#1#2#3{\let\@tempe#1\def\@tempa{#2}\def\@tempb{#3}%
	\futurelet\@tempc\@gifnch}

\def\@gifnch{\ifx\@tempc\@sptoken\let\@tempd\@tempb%
	\else\ifx\@tempc\@tempe\let\@tempd\@tempa\else\let\@tempd\@tempb\fi\fi\@tempd}

\def\SK@set#1{\left\{#1\right\}}
\def\SK@@set#1#2{\{#1\,:\,
	\begin{array}{@{}l@{}}#2\end{array}
	\}}
\def\SK@mset#1{\left\{\!\!\left\{#1\right\}\!\!\right\}}
\def\SK@@mset#1#2{\{\!\!\{#1\,:\,
	\begin{array}{@{}l@{}}#2\end{array}
	\}\!\!\}}
\def\BIG@set#1{\Big\{#1\Big\}}
\def\BIG@@set#1#2{\Big\{#1\:\Big|\:
	\begin{array}{@{}l@{}}#2\end{array}
	\Big\}}
\newcommand{\Set}[1]{\@gifnextchar\bgroup{\SK@@set{#1}}{\SK@set{#1}}}
\newcommand{\Mset}[1]{\@gifnextchar\bgroup{\SK@@mset{#1}}{\SK@mset{#1}}}
\newcommand{\Bigset}[1]{\@gifnextchar\bgroup{\BIG@@set{#1}}{\BIG@set{#1}}}
\makeatother

\title{Small classes have linear neighbourhood complexity}
\author{John Sylvester\thanks{School of Computer Science \& Informatics, University of Liverpool, UK,\\ \texttt{john.sylvester@liverpool.ac.uk}. \orcidlink{0000-0002-6543-2934}
}
	\and
	Viktor Zamaraev\thanks{School of Computer Science \& Informatics, University of Liverpool, UK,\\ \texttt{viktor.zamaraev@liverpool.ac.uk}. \orcidlink{0000-0001-5755-4141}
	}
	\and Maksim Zhukovskii\thanks{School of Computer Science, University of Sheffield, UK, \texttt{m.zhukovskii@sheffield.ac.uk}. \orcidlink{0000-0001-8763-9533}
	}}

\date{}

\begin{document}

\maketitle

\thispagestyle{empty}
    
\begin{abstract}
	A hereditary graph class is \emph{small} if it contains at most $n!c^n$
	labelled graphs on $n$ vertices for some constant $c$. Small classes
    include proper minor-closed classes and classes of bounded twin-width. We prove that every
	small class has \emph{linear neighbourhood complexity}, that is, in every graph
	in the class, vertex neighbourhoods intersect any $k$-vertex set in only
	$O(k)$ distinct ways. This resolves a conjecture of Bonnet, Duron, and the first two authors~(ITCS~2025) and makes the algorithmic, structural, and representation consequences of linear neighbourhood complexity available for every small class. 
	In particular, every small class admits an adjacency labelling scheme
with $O(\log^2 n)$-bit labels, improving the previous $O(\log^3 n)$ bound.
\end{abstract}

\clearpage
\setcounter{page}{1}

\section{Introduction}
\label{sec:introduction}

How much structure is forced by a bound on the number of graphs in a
hereditary class? A graph class is \emph{hereditary} if it is closed under
taking induced subgraphs. For a hereditary class $\mathcal C$, let
$\mathcal C_n$ denote the set of graphs in $\mathcal C$ with vertex set
$[n] := \{1,2, \ldots, n\}$. The function $n\mapsto |\mathcal C_n|$ is called the \emph{speed} of $\mathcal C$.
The classes of complete graphs and edgeless graphs each have speed one,
whereas the class of all graphs has speed $2^{\binom n2}$. Between these
extremes, the speeds of hereditary classes are far from arbitrary.
Classical results show that the low-speed range consists of successive
constant, polynomial, exponential, and factorial layers, with respective
speeds $O(1)$, $n^{\Theta(1)}$, $2^{\Theta(n)}$, and
$2^{\Theta(n\log n)}$~\cite{ScheinermanZito94,Alekseev97,BBW00}.

The layers below factorial admit precise structural descriptions. In
particular, a hereditary class has at most exponential speed if and only
if it has bounded neighbourhood diversity~\cite{Lozin23Parametric}. This
means that the vertex set of every graph in the class can be partitioned
into a bounded number of parts, each inducing a clique or an independent
set, with every two parts either complete or anticomplete to one another. The factorial layer, by contrast, contains many natural classes, including classes of bounded maximum degree, planar graphs, interval
graphs, and disk graphs. Yet no comparable general structural
characterisation of the factorial layer is known.

\paragraph{Small classes.}
A hereditary class $\mathcal C$ is \emph{small} if there exists a constant
$c$ such that $|\mathcal C_n|\leq n!c^n = 2^{(1+o(1))n \log n}$ holds for every $n\geq 1$.
Small classes form a natural extension of the structurally understood
subfactorial regime into the lower part of the factorial layer. They include proper minor-closed classes~\cite{NSTW06}, classes of bounded clique-width~\cite{ALR09}, and, more generally, classes of bounded twin-width~\cite{BGK22}. 

These examples do not capture the full range of small classes. Even
monotone (i.e., classes closed under taking subgraphs) small classes can have unbounded twin-width. More generally, for every countable collection of graph parameters whose boundedness implies smallness, there is a monotone small class on which all these parameters are unbounded~\cite{BDSZZ24}. It is therefore natural to seek useful structural consequences of smallness without expecting them to characterise small classes.

\paragraph{Neighbourhood complexity.}
We focus on the number of distinct neighbourhood patterns that vertices
can realise on a chosen vertex set. 
For a graph $G$ and a set $A\subseteq V(G)$, $\mathcal N_G(A)$ denotes the family of traces on $A$ of all vertex neighbourhoods in $G$, i.e.,
\[
    \mathcal N_G(A)
    :=\{N_G(v)\cap A : v\in V(G)\}.
\]
The \emph{neighbourhood complexity} of $G$ is the function $\nu_G: \mathbb{N} \to \mathbb N$ defined by
\[
    \nu_G(k):=\max\bigl\{|\mathcal N_G(A)|:
        A\subseteq V(G),\ |A|=k\bigr\}.
\]
For a class $\mathcal C$, its neighbourhood complexity is
\[
    \nu_{\mathcal C}(n)
    :=\sup\{\nu_G(n):G\in\mathcal C,\ |V(G)|\geq n\}.
\]
Thus, while speed counts the graphs in a class, neighbourhood complexity
counts the adjacency patterns visible on a set of vertices within those
graphs.

Every class $\cC$ satisfies $\nu_{\mathcal C}(n)\leq 2^n$, with equality for the class of all graphs.
It follows from the Sauer--Shelah lemma \cite{Sau72,She72} that a class $\mathcal C$ has bounded VC-dimension if and only if $\nu_{\mathcal C}(n)=O(n^d)$ for some fixed $d \in \bN$.
We say that $\mathcal C$ has \emph{linear neighbourhood complexity} if
$\nu_{\mathcal C}(n)=O(n)$, where the implicit constant depends only on $\mathcal C$.

Linear neighbourhood complexity is known for classes of bounded
expansion~\cite{RVS19}, bounded twin-width~\cite{BKRTW22,Prz23,BFLP24},
bounded merge-width, and bounded flip-width~\cite{BG25}. Its consequences
include fast graph algorithms~\cite{BDPST26}, short adjacency
labels~\cite{BDSZ25}, and biclique decompositions of small
size~\cite{CMY26}.

\paragraph{Main result.}
A direct connection between smallness and neighbourhood complexity was
established in \cite{BDSZ25}, where it was proved that the neighbourhood complexity of every small class is $O\bigl(n\log n\bigr)$. It was further conjectured that every small class has linear neighbourhood complexity. We prove this conjecture.

\begin{boxmainresult}\label{thm:small-linear-nc}
	Every small graph class has linear neighbourhood complexity.
\end{boxmainresult}

This theorem shows that a bound on the number of graphs forces a uniform
linear bound on the number of neighbourhood patterns on every vertex
set. We note that linear order neighbourhood complexity is {\it best possible}, as witnessed by the small class of graphs of maximum degree at most one.

\paragraph{Complementary results for slow factorial classes.} 

In  light of \cref{thm:small-linear-nc},  it is natural to ask whether linear neighbourhood complexity extends beyond small classes. In particular, we consider hereditary classes of speed $O((n!)^{1+\varepsilon})$, where $\varepsilon>0$ is a fixed constant, which are \emph{slow} factorial but not necessarily small. We show that such classes can have superlinear neighbourhood complexity.

	\begin{restatable}{proposition}{lowerprop}\label{pro:lower}
		 For any fixed $\varepsilon \in (0,1)$  there exists a hereditary class $\mathcal{C}$ with speed $|\mathcal{C}_n| = O((n!)^{1+\varepsilon})$ and neighbourhood complexity  $\Omega( n^{1+\varepsilon^2/64})$.
	\end{restatable}

We also provide an upper bound on the neighbourhood complexity of these classes. 

\begin{restatable}{proposition}{upperprop}\label{pro:upper}
		 For any fixed $\varepsilon \in (0,1)$  and hereditary class $\mathcal{C}$ with speed $|\mathcal{C}_n| \leq  (n!)^{1+\varepsilon} $, the neighbourhood complexity of $\mathcal{C}$ is at most     $ n^{1+ \varepsilon+o(1)}$.
\end{restatable}
    
This result rules out an extreme jump in neighbourhood complexity once one leaves the world of small classes.  Observe also that, as $\varepsilon \to 1$, the bound on neighbourhood complexity approaches $n^2$ as the restriction on speed goes to $(n!)^2$. This aligns with the class of permutation graphs, which have speed $(n!)^{2-o(1)}$ \cite[Theorem 5.2]{BazzaroG09} and neighbourhood complexity $\Omega(n^2)$ \cite[Proposition 14]{FoucaudMNPV17}.

Since the main focus of this paper is linear neighbourhood complexity we did not try to optimise the dependence on $\varepsilon$ in these two propositions. We leave it as an open problem to close the gap between the upper and lower bounds.

\subsection{Applications and discussion}
We combine \cref{thm:small-linear-nc} with known results for classes
of linear neighbourhood complexity to obtain the following corollaries.

\paragraph{Short adjacency labels.}
An \emph{adjacency labelling scheme} for a graph class $\mathcal C$
consists of an encoder and a decoder. Given a graph $G\in\mathcal C$,
the encoder assigns a binary label to each vertex of $G$ so that the
decoder can determine whether two vertices are adjacent from their
labels alone. The decoder is fixed for $\mathcal C$ and does not depend
on the particular graph. The \emph{size} of the scheme is the maximum label
length over all vertices of all $n$-vertex graphs in $\mathcal C$, viewed
as a function of $n$. A class admits an \emph{implicit
representation} if it has an adjacency labelling scheme of size
$O(\log n)$.

A straightforward counting argument shows that every class admitting an implicit
representation has at most factorial speed. The classical \emph{Implicit
Graph Conjecture} posited that this necessary condition is also sufficient
for hereditary graph classes~\cite{Spi03}. Hatami and Hatami refuted the
conjecture by showing that, for every $\delta>0$, there exists a hereditary
factorial class for which every adjacency labelling scheme has size
$\Omega(n^{1/2-\delta})$~\cite{HH22}. 
The Hatami--Hatami counterexample classes are not small, and their approach does not readily yield a small class for which the conjecture fails~\cite{BDSZZ24}.

On the positive side, every monotone small class admits an implicit
representation~\cite{BDSZZ24b}. This motivated the \emph{Small Implicit
Graph Conjecture}, which asserts that every small class admits an
implicit representation~\cite{BDSZZ24b}. The conjecture was subsequently
verified for all \emph{weakly sparse} small classes, namely those
excluding a complete bipartite graph $K_{t,t}$ as a subgraph for some fixed $t$~\cite{BDSZ25}.

Towards this conjecture, it was proved in~\cite{BDSZ25} that every small
class admits an adjacency labelling scheme of size $O(\log^3 n)$.
The proof combines a general method for obtaining short adjacency labels
from low neighbourhood complexity with the structural result that every
small class has neighbourhood complexity $O(n\log n)$. For classes of
linear neighbourhood complexity, the same method gives labels of size
$O(\log^2 n)$~\cite[Theorem~20 and Proposition~6]{BDSZ25}. Our main
theorem therefore provides a further step towards the Small Implicit Graph Conjecture.

\begin{corollary}
\label{th:adjacency-labelling}
    Every small graph class admits an adjacency labelling scheme
    of size $O(\log^2 n)$.
\end{corollary}

\paragraph{Fast algorithms.}
Linear neighbourhood complexity has algorithmic consequences for various problems, including All-Pairs Shortest Paths, matrix multiplication, and triangle detection~\cite{BDPST26}. Our main result makes these algorithmic consequences applicable to all small classes.
\begin{corollary}\label{th:algorithms}
	Let $\cC$ be a small graph class. Then 
	\begin{enumerate}
		\item The All-Pairs Shortest Paths problem can be solved in $O(n^2)$ time for every $n$-vertex graph $G\in\cC$.
		\item The product of the adjacency matrix of any graph $G \in \cC$ with an arbitrary real $n\times n$ matrix can be computed in $O(n^2)$ arithmetic operations\footnote{In the arithmetic model of computation over \(\mathbb R\), where each field operation is assumed to take unit time.}.
		\item For every graph $G\in\cC$ with $n$ vertices and $m$ edges, one can decide whether $G$ contains a triangle by a randomised algorithm in expected time $O(n+m)$.
	\end{enumerate}
\end{corollary}

\paragraph{Bounded symmetric difference.}
Our main theorem also implies that small classes have bounded symmetric difference, i.e., every graph in a small class contains a pair of near-twin vertices.
To state this precisely, for distinct vertices $u,v$ of a
graph $G$, let $\operatorname{sd}_G(u,v)$ be the number of vertices other
than $u$ and $v$ that are adjacent to exactly one of them.
The \emph{symmetric difference} of $G$, denoted $\operatorname{sd}(G)$,
is the minimum $d$ such that every induced subgraph $H$ of $G$ with at
least two vertices contains distinct $u,v\in V(H)$ with
$\operatorname{sd}_H(u,v)\leq d$. 

This parameter was introduced in the study of implicit
representations~\cite{ACLZ15}. It generalises degeneracy, twin-width,
flip-width, and merge-width in the sense that boundedness of any of
these parameters implies bounded symmetric difference.
A simple application of Haussler's packing lemma shows that every class of linear neighbourhood complexity has bounded symmetric difference~(see e.g.~\cite[Appendix~A]{BDPST26}).
Thus, we obtain the following corollary of our main result.

\begin{corollary}\label{th:sym-diff}
	Every small class $\cC$ has bounded symmetric difference. That is, there exists a constant $d:=d(\cC)$ such that $\operatorname{sd}(G)\leq d$ holds for every $G\in\cC$.
\end{corollary}

We note that bounded symmetric difference does not characterise small classes. For example, the class of all graphs of maximum degree at most three has bounded symmetric difference but is not small.

\paragraph{Biclique decompositions of small size.}
A \emph{biclique decomposition} of a graph $G$ is a collection of
bicliques (i.e., complete bipartite subgraphs) whose edge sets partition $E(G)$.
The \emph{size} of a biclique decomposition is the sum of the numbers of vertices in its bicliques. Biclique decompositions of small size have applications in algorithm and data-structure design, extremal combinatorics, and quantum computing (see \cite{CMY26} and references therein). 

It was shown in \cite{CMY26} that every $n$-vertex graph with linear neighbourhood complexity admits a biclique decomposition of size $O(n\log^2 n)$ and such a decomposition can be computed by an efficient randomised algorithm. Together with our main result, this implies the following. 

\begin{corollary}\label{th:biclique-decomp}
    Let $\mathcal C$ be a small graph class. There is an
    $O(n^2\log n)$-time randomised algorithm that, given an $n$-vertex
    graph $G$ in $\mathcal C$, computes a biclique decomposition of $G$ of size
    $O(n\log^2 n)$ with high probability.
\end{corollary}

\paragraph{Distance-$r$ neighbourhood complexity.}
For a positive integer $r$, the \emph{distance-$r$ neighbourhood complexity}
of a graph is defined by replacing open neighbourhoods in the definition of the neighbourhood complexity with \emph{closed $r$-neighbourhoods}
$N_{G}^{(r)}[v]:=\{u\in V(G):\operatorname{dist}_G(u,v)\leq r\}$. 
The \emph{$r$-th power} $G^{(r)}$ of $G$ has vertex set $V(G)$, with distinct
vertices adjacent if and only if their distance in $G$ is at most $r$. 
Open and closed neighbourhood traces on a
set $A$ can differ only for vertices in $A$, so their numbers differ by
at most $|A|$. Consequently, a class has linear distance-$r$ neighbourhood
complexity if and only if the class of its $r$-th powers has linear
neighbourhood complexity.
For monotone graph classes, linear distance-$r$ neighbourhood complexity
for every fixed $r$ characterises bounded expansion~\cite{RVS19}, while
almost-linear distance-$r$ neighbourhood complexity for every fixed $r$
characterises nowhere denseness~\cite{EGKKP16}. Such bounds are an
ingredient in fixed-parameter tractable algorithms for various graph
problems~\cite{EGKKP16,BKRTW22}, including first-order model checking on
monadically stable graph classes~\cite{DreierEMMPT24}.
We leave open whether every small class has linear distance-$r$
neighbourhood complexity for every fixed $r$. This already holds for
weakly sparse small classes, since they have bounded
expansion~\cite{BDSZ25}. A stronger statement that would settle the
question is that smallness is preserved by taking powers and then
hereditary closure. For a hereditary class $\mathcal C$, let
$\mathcal C^{(r)}$ denote the hereditary closure of
$\{G^{(r)}:G\in\mathcal C\}$. Note that taking powers alone cannot increase the number of graphs on any fixed vertex set, since each graph uniquely determines its $r$-th power. The challenge is that the hereditary closure may contain $n$-vertex graphs obtained as induced subgraphs of powers of graphs with more than $n$ vertices.

\begin{problem}
	Is it true that for every small class $\cC$ and every $r \in \bN$, the class $\cC^{(r)}$ is small?
\end{problem}

%%%%%%%%%%%%%%%%%%%%%%%%%%%%%%%%%%
\subsection{Proof strategy of the main result}
%%%%%%%%%%%%%%%%%%%%%%%%%%%%%%%%%%
Our approach extends the counting argument used in~\cite[Theorem~4]{BDSZ25} to prove the $O(n\log n)$ bound on the neighbourhood complexity of small classes. We first outline this argument, then describe an extension that yields a linear bound on neighbourhood complexity in a special case. Finally, we explain the strategy for the general case.

\paragraph{Basic approach.}
For a bipartite graph with parts $A$ and $B$, we say that $B$ is \emph{$A$-diverse} if its vertices have pairwise distinct neighbourhoods in~$A$.

Let $\cC$ be a small class with $|\cC_m|\leq m!c^m$ for every $m$, where $c\geq 1$. Assume that the neighbourhood complexity of $\cC$ is not $O(n \log n)$. Then, for any $C > 0$, there exist a graph $G\in\cC$ and an $n$-vertex set $A$ in $G$ such that the vertices of $G$ have at least $2CN$ distinct neighbourhoods in $A$, where $N = n\log n$. In particular, there is a set $B$ disjoint from $A$ such that $|B|=CN$ and $B$ is $A$-diverse. To obtain a contradiction, we show that $G$ gives rise to more than $N!c^N$ labelled $N$-vertex graphs in $\cC$. To do so,
we fix a labelling of $A$ by $[n]$, select a subset $\tilde B$ of $B$ of size $N-n$, and label its vertices by $\{n+1,\ldots,N\}$ (see \cref{fig:proof-strategy-basic}). Since $B$ is $A$-diverse and the labelling of $A$ is fixed, different choices of $\tilde B$ and its labelling give pairwise distinct labelled bipartite graphs $G[A,\tilde B]$ spanned by the edges of $G$ between $A$ and $\tilde B$. Restoring the edges of $G$ within each part preserves these differences and yields graphs in $\cC_N$. Thus we obtain
\[
\binom{CN}{N-n}(N-n)!>N!C^{N/3}>N!c^N
\]
distinct graphs in $\cC_N$, contradicting the smallness assumption when $C > c^3$.

% Proof-strategy diagrams: a standalone basic case and two extension panels.
% This file defines shared styles and one command for each diagram.

% A core colour used only in the proof-strategy overview. Later figures retain the
% paper's established indigo/blue/peach set palette.
\definecolor{StrategyCore}{RGB}{88,143,118}

\tikzset{
	strategy/ambient/.style={
		draw=black!38,
		fill=black!2,
		line width=.22pt,
		rounded corners=4pt,
		outer sep=0pt,
		inner sep=0pt
	},
	strategy/ambient faded/.style={
		draw=black!17,
		fill=black!1,
		line width=.18pt,
		rounded corners=4pt,
		outer sep=0pt,
		inner sep=0pt
	},
	strategy/pool/.style={
		strategy/ambient
	},
	strategy/extension/.style={
		draw=black!60,
		fill=black!2,
		line width=.95pt,
		rounded corners=3pt,
		outer sep=0pt,
		inner sep=0pt
	},
	strategy/core/.style={
		strategy/ambient,
		fill=none
	},
	strategy/core pool/.style={
		strategy/pool,
		fill=none
	},
	strategy/core link/.style={
		fill=StrategyCore,
		fill opacity=.22,
		draw=none
	},
	strategy/selected/.style={
		draw=black!60,
		line width=.95pt
	},
	strategy/label/.style={
		font=\fontsize{8pt}{9.2pt}\selectfont,
		anchor=north west,
		inner sep=0pt
	}
}

% Caption key matching the fill, border colour, and weight of \tilde B.
\DeclareRobustCommand{\proofstrategyextensioncaptionbox}[1]{%
	\tcbox[on line, colback=black!2, colframe=black!60, boxrule=.95pt,
		arc=3pt, boxsep=2pt, left=2pt, right=2pt,
		top=.3pt, bottom=.3pt]{#1}%
}

% Wrap both core sides in a straight-sided ribbon with a 1mm margin.
% The set outlines are restored afterward to keep them clear.
\newcommand{\proofstrategycorelink}[2]{%
	\path[strategy/core link, rounded corners=3.5pt]
		($(#1.north west)+(-1mm,1mm)$)
		-- ($(#1.north east)+(1mm,1mm)$)
		-- ($(#1.south east)+(1mm,-1mm)$)
		-- ($(#2.north east)+(1mm,1mm)$)
		-- ($(#2.south east)+(1mm,-1mm)$)
		-- ($(#2.south west)+(-1mm,-1mm)$)
		-- ($(#2.north west)+(-1mm,1mm)$)
		-- ($(#1.south west)+(-1mm,-1mm)$)
		-- cycle;
}

\newcommand{\proofstrategyredraw}[2]{%
	\path[#1] (#2.south west) rectangle (#2.north east);
}

% Basic O(n log n) count, shown as a standalone figure.
\newcommand{\proofstrategybasic}{%
	\begin{tikzpicture}
		\path[use as bounding box] (-.65cm,-1.55cm) rectangle (8.85cm,1.35cm);

		\node[
			strategy/ambient,
			minimum width=8.70cm,
			minimum height=12.5mm
		] (strategy-a-B) at (4.25cm,6mm) {};

		\node[
			strategy/extension,
			minimum width=6.40cm,
			minimum height=8.5mm
		] (strategy-a-Btilde) at (4.75cm,6mm) {};

		\node[
			strategy/ambient,
			strategy/selected,
			minimum width=3.60cm,
			minimum height=10.5mm
		] (strategy-a-A) at (4.25cm,-9.5mm) {};

		\node[font=\small, anchor=east]
			at ($(strategy-a-B.west)+(-2.5mm,0)$) {$B$};
		\node[font=\small, anchor=east]
			at ($(strategy-a-A.west)+(-2.5mm,0)$) {$A$};
		\node[strategy/label]
			at ($(strategy-a-Btilde.north west)+(1.7mm,-1.5mm)$) {$\tilde B$};
	\end{tikzpicture}%
}

% (a) Bounded A-neighbourhoods. The ribbon groups A and B'' into the core
% F[A,B'']; A and \tilde B have heavy boundaries because together they
% induce the selected extension.
\newcommand{\proofstrategybounded}{%
	\begin{tikzpicture}
		\path[use as bounding box] (.45cm,-1.65cm) rectangle (7.60cm,1.35cm);

		\node[
			strategy/ambient,
			minimum width=6.40cm,
			minimum height=12.5mm
		] (strategy-b-B) at (4.25cm,6mm) {};

		\node[
			strategy/extension,
			minimum width=4.70cm,
			minimum height=8.5mm
		] (strategy-b-Btilde) at (4.85cm,6mm) {};

		\node[
			strategy/core,
			minimum width=1.35cm,
			minimum height=5.7mm
		] (strategy-b-Bdoubleprime) at (6.30cm,5.5mm) {};

		\node[
			strategy/core,
			strategy/selected,
			minimum width=3.30cm,
			minimum height=10.5mm
		] (strategy-b-A) at (4.25cm,-9.5mm) {};

		\proofstrategycorelink{strategy-b-Bdoubleprime}{strategy-b-A}
		\proofstrategyredraw{strategy/ambient, fill=none}{strategy-b-B}
		\proofstrategyredraw{strategy/extension, fill=none}{strategy-b-Btilde}
		\proofstrategyredraw{strategy/core}{strategy-b-Bdoubleprime}
		\proofstrategyredraw{strategy/core, strategy/selected}{strategy-b-A}

		\node[font=\small, anchor=east]
			at ($(strategy-b-B.west)+(-2.5mm,0)$) {$B$};
		\node[font=\small, anchor=east]
			at ($(strategy-b-A.west)+(-2.5mm,0)$) {$A$};
		\node[strategy/label]
			at ($(strategy-b-Btilde.north west)+(1.7mm,-1.5mm)$) {$\tilde B$};
		\node[font=\fontsize{8pt}{9.2pt}\selectfont, inner sep=0pt]
			at (strategy-b-Bdoubleprime.center) {$B''$};
	\end{tikzpicture}%
}

% (b) General case. The outer A and B are merely ambient supersets and are
% deliberately faded; B' remains visible as the pool for the chosen \tilde B.
\newcommand{\proofstrategygeneral}{%
	\begin{tikzpicture}
		\path[use as bounding box] (.45cm,-1.65cm) rectangle (7.60cm,1.35cm);

		\node[
			strategy/ambient faded,
			minimum width=6.40cm,
			minimum height=14mm
		] (strategy-c-B) at (4.25cm,6mm) {};

		\node[
			strategy/pool,
			minimum width=5.60cm,
			minimum height=11.2mm
		] (strategy-c-Bprime) at (4.45cm,6mm) {};

		\node[
			strategy/extension,
			minimum width=4.40cm,
			minimum height=8.5mm
		] (strategy-c-Btilde) at (4.90cm,6mm) {};

		\node[
			strategy/core pool,
			minimum width=2.20cm,
			minimum height=5.7mm
		] (strategy-c-Bdoubleprime) at (5.85cm,5.5mm) {};

		\node[
			strategy/ambient faded,
			minimum width=4.00cm,
			minimum height=12mm
		] (strategy-c-A) at (4.25cm,-9.5mm) {};

		\node[
			strategy/core,
			strategy/selected,
			minimum width=2.20cm,
			minimum height=8.5mm
		] (strategy-c-Aprime) at (4.55cm,-9.5mm) {};

		\proofstrategycorelink{strategy-c-Bdoubleprime}{strategy-c-Aprime}
		\proofstrategyredraw{strategy/ambient faded, fill=none}{strategy-c-B}
		\proofstrategyredraw{strategy/pool, fill=none}{strategy-c-Bprime}
		\proofstrategyredraw{strategy/extension, fill=none}{strategy-c-Btilde}
		\proofstrategyredraw{strategy/ambient faded, fill=none}{strategy-c-A}
		\proofstrategyredraw{strategy/core pool}{strategy-c-Bdoubleprime}
		\proofstrategyredraw{strategy/core, strategy/selected}{strategy-c-Aprime}

		\node[font=\small, text=black!42, anchor=east]
			at ($(strategy-c-B.west)+(-2.5mm,0)$) {$B$};
		\node[font=\small, text=black!42, anchor=east]
			at ($(strategy-c-A.west)+(-2.5mm,0)$) {$A$};
		\node[strategy/label]
			at ($(strategy-c-Bprime.north west)+(1.5mm,-1.3mm)$) {$B'$};
		\node[strategy/label]
			at ($(strategy-c-Btilde.north west)+(1.7mm,-1.5mm)$) {$\tilde B$};
		\node[font=\fontsize{8pt}{9.2pt}\selectfont, inner sep=0pt]
			at (strategy-c-Bdoubleprime.center) {$B''$};
		\node[strategy/label]
			at ($(strategy-c-Aprime.north west)+(1.7mm,-1.5mm)$) {$A'$};
	\end{tikzpicture}%
}

\begin{figure}[htbp]
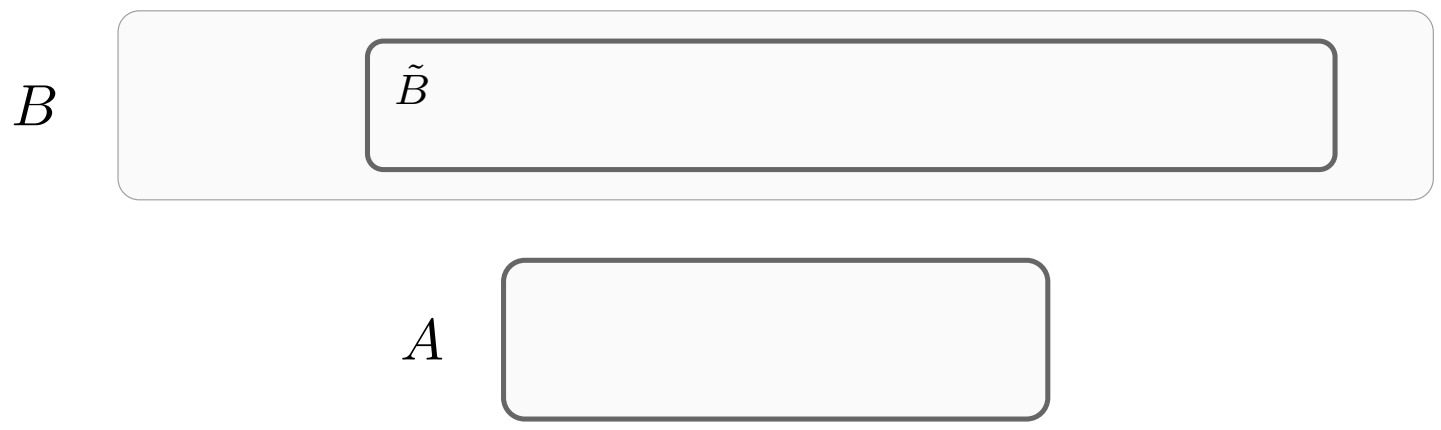

	\centering
	\proofstrategybasic
	\caption{The basic approach. We count labelled subgraphs $G[A,\tilde B]$ by choosing and labelling $\tilde B\subseteq B$, while keeping the labelling of $A$ fixed.} 
	\label{fig:proof-strategy-basic}
\end{figure}

In this form, the approach cannot improve the order of the bound on neighbourhood complexity $f(n) = |B|$. With the labelling of $A$ fixed, selecting and labelling
$k$ vertices of $B$ gives $\binom{f(n)}{k}k!$ choices.
To contradict smallness using this count, we would need
\[
    \binom{f(n)}{k}k!>(n+k)!c^{n+k},
\]
which can hold only if $f(n)=\Omega(n\log n)$.

\paragraph{Extension of the basic approach: bounded-degree case.}
A natural idea for obtaining a linear bound is to vary the labellings of $A$. However, different choices of a subset $\tilde B\subseteq B$ and a labelling of $A\cup\tilde B$ may produce the same labelled graph. The challenge is to control this multiplicity when both $|B|$ and $|\tilde B|$ are linear in $n$.

\begin{figure}[htbp]
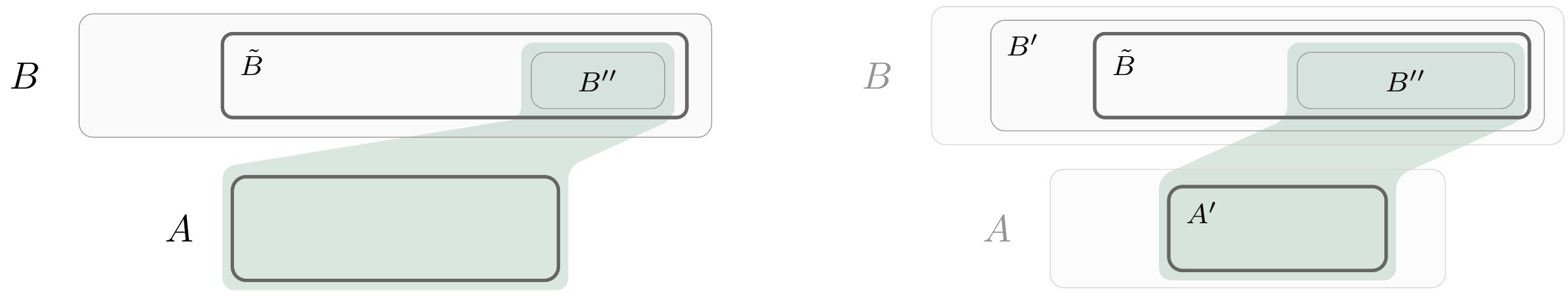

	\centering
	\begin{subfigure}[t]{.48\linewidth}
		\centering
		\proofstrategybounded
		\caption{\footnotesize Illustrative special case: vertices in $B$ have bounded degree. We find a core subgraph \setcaptionbox{StrategyCore!22}{F[A,B'']} in $F$, then count its labelled extensions $F[A,\tilde B]$ with $B''\subseteq\tilde B\subseteq B$.}
		\label{fig:proof-strategy-bounded}
	\end{subfigure}\hfill%
	\begin{subfigure}[t]{.48\linewidth}
		\centering
		\proofstrategygeneral
		\caption{\footnotesize General case. We first find a core subgraph \setcaptionbox{StrategyCore!22}{F[A',B'']} inside a suitable subgraph $F[A',B']$ of $F$, then count its labelled extensions $F[A',\tilde B]$ with $B''\subseteq\tilde B\subseteq B'$.}
		\label{fig:proof-strategy-general}
	\end{subfigure}
	\caption{Extensions of the basic approach. We find a \setcaptionbox{StrategyCore!22}{\text{core subgraph}} with few automorphisms and count its labelled extensions. The automorphism bound controls how many choices of an extension and a labelling produce the same labelled graph, giving a lower bound on the number of distinct labelled graphs.}
	\label{fig:proof-strategy}
\end{figure}

Our key observation is that an induced subgraph with few automorphisms can control this multiplicity. 
We illustrate this idea in the special case where vertices in $B$ have bounded degree\footnote{We discuss the bounded-degree case here only to illustrate the counting strategy; the full proof does not treat it separately.}. More specifically, consider $F:=G[A,B]$ for a graph $G\in\cC$, where $A,B$ are disjoint, $|A|=n$, $|B|=Cn$, $C$ is sufficiently large, and $B$ is $A$-diverse. Under suitable minimality assumptions, \cref{lm:bounded} gives a set $B''\subseteq B$ of size less than $n$ such that $F[A,B'']$ has at most $e^{O(n\ln C)}$ automorphisms. We call $F[A,B'']$ the \emph{core} (see \cref{fig:proof-strategy}\textup{(a)}). 

We count labelled \emph{extensions} $F[A,\tilde B]$ of the core such that $B''\subseteq\tilde B\subseteq B$, $|\tilde B| = \sqrt{C}n$, $A$ is labelled by $[n]$, and $\tilde B$ is labelled by $\{n+1, \ldots, \sqrt{C}n + n\}$. To do this, fix one extension and one such labelling. The choices of an extension and a labelling that produce this same labelled graph correspond to isomorphisms
from $F[A,\tilde B]$ onto induced subgraphs of $F$ that contain $B''$, with each isomorphism preserving $A$ setwise. There are at most $e^{O(n\ln C)}$ such isomorphisms:
\begin{itemize}
	\item there are at most $\binom{|B|}{|B''|}=e^{O(n\ln C)}$ choices for the image of $B''$; 
	\item once this image is fixed, there are at most $|\mathrm{Aut}(F[A,B''])|=e^{O(n\ln C)}$ choices for the restriction of the isomorphism to $A\cup B''$;
	\item this restriction determines the images of the remaining vertices, since $B$ is $A$-diverse.
\end{itemize}
On the other hand, since $|B''|<n$, the number of extensions with $|\tilde B|=\sqrt{C}n$ is at least
\[
\binom{(C-1)n}{(\sqrt{C}-1)n}=e^{\Omega(\sqrt{C}n\ln C)}.
\]
Each extension admits $n!(\sqrt{C}n)!$ labellings with the prescribed label ranges. Setting $N:=(\sqrt{C}+1)n$, we have $n!(\sqrt{C}n)!\geq N!/2^N$. Thus the number of extensions absorbs both the factor $2^N$ and the multiplicity $e^{O(n\ln C)}$, leaving more than $N!c^N$ distinct labelled graphs for sufficiently large $C$. As before, this contradicts smallness after restoring the edges within each part.

\paragraph{Extension of the basic approach: general case.}
In the general case, where the degrees of vertices in $B$ are not bounded, we find the core $F[A',B'']$ inside a suitable subgraph $F[A',B']$ of $F$ (see \cref{fig:proof-strategy}\textup{(b)}). More specifically, our main technical result, \cref{thm:auto}, gives sets $A'\subseteq A$ and $B''\subseteq B'\subseteq B$ such that $F[A',B'']$ has few automorphisms and $B'\setminus B''$ is $A'$-diverse. It also gives the bounds on the sizes of these sets needed for the counting argument. We again count labelled extensions $F[A',\tilde B]$ with $B''\subseteq\tilde B\subseteq B'$. Although diversity is guaranteed only outside the core, requiring every extension to contain $B''$ still ensures that each isomorphism considered in the count is determined by its restriction to the core. The resulting multiplicity bound yields the contradiction to smallness (see \cref{sec:MainTheoremFromLNC} for details).

It remains to explain how we find a suitable core in the general case. We first choose an inclusion-maximal set $R\subseteq B$ whose vertices have neighbourhoods in $A$ with large pairwise symmetric differences. The packing lemma (\cref{thm:packing}) bounds the size of $R$. We partition $A$ according to neighbourhoods in $R$ and choose one representative from each class. The \emph{canonical neighbourhood} of a vertex $v\in B\setminus R$ is the union of the classes whose representatives are adjacent to $v$. We show that a single choice of representatives ensures that the symmetric difference between the actual and canonical neighbourhoods has bounded size for a constant fraction of the vertices in $B$. We call these vertices \emph{good} and denote the set of good vertices by $B_{\good}$. The \emph{profile} of a good vertex is its neighbourhood in the set of representatives, and its \emph{defect} is the symmetric difference between its actual and canonical neighbourhoods.

The proof distinguishes two cases according to the number of distinct defects among good vertices.
When there are many distinct defects, we build an auxiliary bipartite graph with parts $A$ and $B_{\good}$ in which each good vertex is adjacent to every vertex in its defect. In this graph, every vertex in $B_{\good}$ has a neighbourhood of bounded size, allowing us to apply \cref{lm:bounded}. After a suitable minimality restriction, this lemma provides a core subgraph in the auxiliary graph with few automorphisms. 

When there are few distinct defects, we instead build an auxiliary bipartite graph whose vertices are profiles and defects. A profile vertex is adjacent to a defect vertex if there exists a good vertex with that (profile, defect) pair. 
This graph has many edges relative to its number of vertices. Repeatedly deleting vertices of small degree leaves a non-empty induced subgraph of large minimum degree. Applying \cref{lm:tree} then gives an induced subgraph with a spanning tree of bounded degree. By \eqref{eq:tree-automorphisms}, this tree has few automorphisms and serves as a core subgraph of the auxiliary graph.

In both cases, we use the auxiliary bipartite graph to define $A',B',B''$ and bound the automorphism group of the core subgraph $F[A',B'']$ of the original graph. To obtain this bound, we select a bounded number of vertex sets in $A'\cup B''$ whose common setwise stabiliser $\Gamma$ in $\mathrm{Aut}(F[A',B''])$ admits an injective homomorphism into the automorphism group of the core of the auxiliary graph. By the orbit--stabiliser theorem, $\mathrm{Aut}(F[A',B''])$ and $\Gamma$ differ in size by at most a factor exponential in $|A'|+|B''|$ (\cref{lm:orbit-stabiliser}). Together with the bound on the automorphism group of the core of the auxiliary graph, this gives the desired bound on $|\mathrm{Aut}(F[A',B''])|$.

%%%%%%%%%%%%%%%%%%%%%%%%%%%%%%%%%%
\section{Preliminaries}
%%%%%%%%%%%%%%%%%%%%%%%%%%%%%%%%%%

For two sets $A,B$, let $A\triangle B = (A\setminus B)\cup (B\setminus A) $ denote their symmetric difference. 

For a graph $G$, we denote by $V(G)$ and $E(G)$ the set of vertices and the set of edges of $G$, respectively.
For a vertex $v$ in $V(G)$, we denote by $N(v)$ the neighbourhood of $v$. For a set $A \subseteq V(G)$, the \emph{$A$-neighbourhood of $v$ in $G$} is the set $N_A(v) := N(v) \cap A$. When $G$ is clear from the context, for the sake of brevity, we  write $A$-neighbourhood. 

For a graph $G$ and a set $U\subseteq V(G)$, we denote by $G[U]$ the induced subgraph of $G$ on the set of vertices $U$.
If $A$ and $B$ are disjoint subsets of $V(G)$, we denote by $G[A,B]$ the bipartite graph with parts $A$ and $B$ such that $a\in A$ and $b\in B$ are adjacent in $G[A,B]$ if and only if they are adjacent in $G$.

Recall, for a hereditary class $\cC$, we denote by
$\cC_n$ the set of graphs in $\mathcal C$ with vertex set $[n] $. That is, graphs in $\cC_n$ are labelled and two isomorphic graphs in $\cC_n$ are considered to be different if they have different labellings. 
A class $\cC$ is \emph{small} if there exists a constant $c$ such that $|\cC_n| \leq n!c^n$ holds for all $n \in \bN$.

%%%%%%%%%%%%%%%%%%%%%%%%%%%%%%%%%%
\subsection{Packing lemma}
%%%%%%%%%%%%%%%%%%%%%%%%%%%%%%%%%%
Let $A$ be an $n$-element set and let $\mathcal H\subseteq 2^A$ be a family of its subsets. For a positive integer $C$, we say that $\mathcal{H}$ is {\it $C$-separated} if $|V\triangle V'|>C$ for every pair of different $V,V'\in\mathcal{H}$. The {\it shatter function} of $\mathcal{H}$ is the function $\tau_{\mathcal{H}}:\{0,1,\ldots,n\}\to\mathbb{N}$ defined by 
$$
 \tau_{\mathcal{H}}(k)=\max\{|\{V\cap U:\,V\in\mathcal{H}\}|:\,U\subseteq A, |U|=k\}.
$$

In particular, the neighbourhood complexity of $G$ is the shatter function of the set system consisting  of neighbourhoods in $G$. The \emph{Vapnik–Chervonenkis dimension} $d$ (or \emph{VC dimension} for short) of $\mathcal{H}$ is the maximum cardinality  of a set $U\subseteq A$ such that \[\{V\cap U:\,V\in\mathcal{H}\}=2^U.\]  We will use a corollary of the following packing lemma~\cite[Lemma 23]{BDSZ25}, which is a shatter-function formulation of Haussler's packing lemma~\cite[Theorem 1]{Haussler}.

\begin{theorem}[Packing lemma]
\label{thm:packing}
Let $1\leq C\leq n$ be integers. Let $A$ be a set of size $n$ and let $\mathcal{H}$ be a non-empty family of subsets of $A$ of VC-dimension $d\geq 1$. If $\mathcal{R}\subseteq\mathcal{H}$ is $C$-separated, then for $k:=\lceil4dn/C\rceil \leq n $ we have $|\mathcal{R}|\leq 2\tau_{\mathcal{H}}(k)$.
\end{theorem}

This gives the following result in our setting.  
\begin{corollary}
\label{cor:packing}
Let $C$ be a large integer. Let $A$ be a set of size $n>C^{10}$. Let $\mathcal{H}$ be a family of subsets of $A$ such that, for every non-empty $U\subseteq A$, $|\{V\cap U:\,V\in\mathcal{H}\}|\leq C|U|$. If $\mathcal{R}\subseteq \mathcal{H}$ is $C^{10}$-separated, then $|\mathcal{R}|\leq\frac{n}{C^7}$.
\end{corollary}

\begin{proof}
Let $k=\lceil 4dn/C^{10}\rceil$, where $d$ is the VC-dimension of $\mathcal{H}$. If $d=0$, then $|\mathcal{R}|\leq 1$, so the claim is immediate. Assume $d\geq 1$.  Let $U\subseteq A$ be a set of size $d$ that is shattered by $\mathcal{H}$. Since $2^d=|\{V\cap U:\,V\in\mathcal{H}\}|\leq Cd$ and $C$ is large, we get that $4d/C^{10}\leq 2/C^9<1$, implying $k\leq n$. By \cref{thm:packing}, we have
$$
|\mathcal{R}|\leq 2\tau_{\mathcal{H}}(k)\leq 2Ck\leq 8dn/C^9+2C\leq 4n/C^8+2n/C^9\leq n/C^7,
$$
as desired.
\end{proof}

\subsection{Diverse sets}

We recall that for sets $A,B \subseteq V(G)$ in a graph $G$, we say that $B$ is \emph{$A$-diverse in $G$} if vertices in $B$ have pairwise distinct $A$-neighbourhoods in $G$, that is all sets $N_A(b), b \in B$ are pairwise different. When $G$ is clear from the context we  write $A$-diverse for brevity. We say that an $A$-diverse set $B$, where $A\neq \emptyset$, is \emph{minimally $A$-diverse} if, for every nonempty proper subset $U\subset A$,
\[
|\{N_U(b):b\in B\}|<
\frac{|B|}{|A|}|U|.
\]

\begin{lemma}
\label{lm:small-diverse}
Let $F$ be a bipartite graph with non-empty parts $A$ and $B$, and let $B$ be $A$-diverse. Then there exists a set $A'\subseteq A$ with $|A'|\leq |B|-1$ such that $B$ is $A'$-diverse and  $A'$ is $B$-diverse. 
\end{lemma}

\begin{proof}
We construct the set $A'$ iteratively. Start from an empty set. Then we add to $A'$ one element at a time, so that this element refines the partition of $B$ with respect to its adjacency relations with $A'$. Unless $B$ is completely refined, at each step we can find such an element since $B$ is $A$-diverse. The total number of refinement steps is at most $|B|-1$, so $|A'|\leq|B|-1$. To conclude, let us show that $A'$ is $B$-diverse. Assume $u,v\in A'$ have equal neighbourhoods in $B$ and $v$ was added later than $u$. Then $v$ does not refine the partition of $B$ obtained just before the step when $v$ was added --- a contradiction.
\end{proof}

\begin{lemma}
\label{lm:minimal-diverse}
Let $C$ be a positive integer, and let $F$ be a bipartite graph with non-empty parts $A$ and $B$ such that $B$ is $A$-diverse and $|B|\geq C|A|$. Then there exist non-empty sets $A'\subseteq A$ and $B'\subseteq B$ such that $B'$ is minimally $A'$-diverse and $|B'|=C|A'|$.
\end{lemma}

\begin{proof}
The proof is straightforward: find a non-empty inclusion-minimal subset $A'\subseteq A$ such that $B$ has at least $C|A'|$ traces on $A'$. Then retain exactly $C|A'|$ vertices in $B$ such that their neighbourhoods in $A'$ are different.
\end{proof}

%%%%%%%%%%%%%%%%%%%%%%%%%%%%%%%%%%
\subsection{Groups and automorphisms}
\label{sec:groups}
%%%%%%%%%%%%%%%%%%%%%%%%%%%%%%%%%%
An \emph{action} of a finite group $H$ on a finite set $\Gamma$ assigns to each $h\in H$ a permutation $x\mapsto hx$ of $\Gamma$, such that $1_Hx=x$ and $(gh)x=g(hx)$ for all $g,h\in H$ and $x\in\Gamma$. Here $1_H$ is the identity element of $H$. The \emph{order} of $H$ is the number of its elements $|H|$. We write $\mathrm{Sym}(\Gamma)$ for the group of all permutations of $\Gamma$ under composition; its subgroups are \emph{permutation groups} on $\Gamma$. For a graph $F$, its \emph{automorphism group} $\mathrm{Aut}(F)$ consists of the permutations of $V(F)$ that preserve adjacency. For a directed graph, automorphisms also preserve the directions of arcs.

For $S\subseteq\Gamma$ and $h\in H$, write $hS:=\{hx:\,x\in S\}$. Given sets $S_1,\ldots,S_k\subseteq\Gamma$, their \emph{setwise stabiliser} is a subgroup
\[
    H_{S_1,\ldots,S_k}:=\{h\in H:\,hS_i=S_i\text{ for every }i\in[k]\}.
\]
The \emph{pointwise stabiliser} of $S\subseteq\Gamma$ is a subgroup
\[
    H_{(S)}:=\{h\in H:\,hx=x\text{ for every }x\in S\}.
\]
The \emph{orbit} of the ordered tuple $(S_1,\ldots,S_k)$ is
\[
    \mathrm{Orb}_H(S_1,\ldots,S_k):=\{(hS_1,\ldots,hS_k):\,h\in H\}.
\]

A \emph{group homomorphism} $\varphi:H\to K$ is a map satisfying $\varphi(gh)=\varphi(g)\varphi(h)$ for all $g,h\in H$. Its \emph{kernel} is $\ker\varphi:=\{h\in H:\,\varphi(h)=1_K\}$. We recall the standard fact that $\varphi$ is injective if and only if its kernel is trivial, that is, $\ker\varphi=\{1_H\}$. In this case, $|H|\leq |K|$.

We will use the following form of the orbit--stabiliser theorem to bound the order of a group by the order of a subgroup preserving several specified sets.
\begin{lemma}[Orbit--stabiliser bound]
\label{lm:orbit-stabiliser}
Let a finite group $H$ act on a finite set $\Gamma$, and let $S_1,\ldots,S_k\subseteq\Gamma$. Then
$
|H|=|\mathrm{Orb}_H(S_1,\ldots,S_k)|\cdot|H_{S_1,\ldots,S_k}|
\leq (2^k)^{|\Gamma|}\cdot|H_{S_1,\ldots,S_k}|.
$
\end{lemma}

The next lemma provides a small core with few automorphisms when neighbourhoods on one side have bounded size. It will be a key ingredient in the proof of \cref{thm:auto}.

\begin{lemma}
\label{lm:bounded}
Let $C,C'$ be large integers. Let $F$ be a bipartite graph with parts $A$ and $B$, where $B$ is minimally $A$-diverse and $|B|=C|A|$. Let every neighbourhood $N_A(v)$, $v\in B$, have size at most $C'$. Then there exists a set $B'\subseteq B$ of size at most $|A|-1$ such that $|\mathrm{Aut}(F[A, B'])|\leq \max\{C,C'\}^{4|A|}.$ 
\end{lemma}

\begin{proof}
Consider a blow-up $F^{\star}$ of $F$ where every vertex of $A$ is replaced by $C$ copies with the same neighbours in $B$. Let $A^{\star}$ be the corresponding blow-up of $A$. Let us show that there exists a perfect matching from $B$ to $A^{\star}$ in $F^{\star}$, unless $B$ contains an isolated vertex. If the latter is indeed the case, we show the existence of a matching that saturates $B\setminus\{b_0\}$, where $b_0$ is the (unique) isolated vertex in $B$. In general let 
\[Z =\begin{cases}
\emptyset & \text{if $B$ has no isolated vertex,}\\ 
\{b_0\} & \text{otherwise}
\end{cases}.\] It suffices to verify Hall's condition. Take $V\subseteq B\setminus Z$ and assume that its neighbourhood $U^{\star}$ in $A^{\star}$ has size less than $|V|$. Let $U\subseteq A$ be the set of original vertices of $U^{\star}$. 
Since $V\subseteq B\setminus Z$, the set $U$ is nonempty. Moreover, $C|U|<|V|\le C|A|$, so $U\subset A$ is proper. As the vertices of $V$ have distinct neighbourhoods contained in $U$, there are at least  $C|U|$ distinct traces on $U$, contradicting minimal $A$-diversity. Therefore, $F^{\star}$ contains a matching $P$ saturating $B\setminus Z$. Consider now a spanning oriented subgraph $F^{\downarrow}\subseteq F$, where $(a,b)\in F^{\downarrow}$ if $\{a^{\star},b\}\in P$ for some copy $a^{\star}$ of $a$ in $A^{\star}$. Every vertex of $B\setminus Z$ has indegree one in $F^{\downarrow}$, while $b_0$, if it exists, has indegree zero. If $Z=\emptyset$, every vertex of $A$ has outdegree $C$. Otherwise, exactly one vertex of $A$ has outdegree $C-1$, while all other vertices of $A$ have outdegree $C$.

We construct the set $B'$ vertex by vertex. 
Choose $b_1 \in B$ with $|N_A(b_1)| \geq 2$ and let $U_1 = N_A(b_1)$; such a vertex exists since $B$ is $A$-diverse and $C|A|=|B|>|A|+1$. For $i \geq 2$, find a vertex $b_i$ such that the set of neighbours of $b_i$ in $F$ is not entirely inside the set $U_{i-1}$ of neighbours of $\{b_1,\ldots,b_{i-1}\}$ in $F$ whereas the only in-neighbour of $b_i$ in $F^{\downarrow}$ lies in $U_{i-1}$. Observe that this process halts once $N_A(\{b_1,\ldots,b_i\})=A$. Indeed, otherwise, there exists a proper subset $U\subset A$ such that every vertex in $B$ with its $F^{\downarrow}$-in-neighbour in $U$ has all its $F$-neighbours in $U$. The set $U$ has either exactly $C|U|$ out-neighbours in $B$ if $B$ does not have isolated vertices or at least $C|U|-1$ out-neighbours in $B$ if $b_0\in B$. If every such out-neighbour has its entire $F$-neighbourhood in $U$, this contradicts minimality of $A$. In particular, in the case of $C|U|-1$ out-neighbours, the number of vertices in $B$ with different neighbourhoods in $U$ is $C|U|$ as well, since $b_0$ is among these vertices and contributes the empty neighbourhood. We get $|B'|\leq|A|-1$. Observe $b_0\notin B'$. Let $F^{\uparrow}$ be the induced subgraph $F[A, B']$, where edges are oriented from $B'$ to $A$. 

Note that for every vertex $v\in A\cup B'$, there exists a directed path from $b_1$ to $v$ in $F^{\updownarrow}: =F^{\uparrow}\cup F^{\downarrow}[A, B']$. This follows by induction since from $\{b_1,\ldots,b_{i-1}\}$ every vertex in $U_{i-1}$ is reachable, and there exists a vertex in $U_{i-1}$, from which $b_i$ is reachable. Since the maximum outdegree of $F^{\updownarrow}$ is at most $\Delta:=\max\{C,C'\}$, we get $|\mathrm{Aut}(F^{\updownarrow})|\leq 2|A|\cdot \Delta^{2|A|}$. Indeed, there are at most $|B'|+|A|\leq 2|A|$ ways to choose the image of $b_1$. Then, fixing a spanning arborescence of $F^{\updownarrow}$ rooted at $b_1$ and processing vertices in top-down order, each vertex must map to an out-neighbour of its parent’s image, giving at most $\Delta$ choices per vertex.

We note that $F':=F[A, B']$ is the undirected version of $F^{\updownarrow}$, which (by the above) is connected. Since $|B'|<|A|$, its two bipartition classes have different sizes. Consequently, every automorphism of $F'$ preserves $A$ and $B'$ setwise. Since every automorphism of $F'$ preserves the bipartition and since every edge $(a,b)$ in $F^{\downarrow}[A,  B']$ has its copy $(b,a)$ directed in the opposite way in $F^{\uparrow}$, an automorphism of $F'$ preserves the edges of $F^{\downarrow}[A, B']$ if and only if it is an automorphism of $F^{\updownarrow}$.

The number of ways to choose a single edge from every vertex of $B'$ to $A$ is at most $(C')^{|B'|}<\Delta^{|A|}$. Therefore, the orbit of the set of (undirected) edges of $F^{\downarrow}[A, B']$ in $\operatorname{Aut}(F')$ has size at most $\Delta^{|A|}$. By the orbit--stabiliser theorem, \[|\mathrm{Aut}(F')|<2|A|\cdot \Delta^{2|A|}\cdot \Delta^{|A|}<\Delta^{4|A|} = (\max\{C,C'\})^{4|A|},\] completing the proof. 
\end{proof}

Our second source of graphs with few automorphisms is trees of bounded degree. The following lemma\footnote{\cite[Lemma 12]{BDSZZ24b} is stated for a graph $G$ of minimum degree \emph{exactly} $d$. However, the version stated here follows by first passing to a subgraph of minimum degree exactly $d$, applying \cite[Lemma 12]{BDSZZ24b}, and then adding any missing edges back on at the end.} supplies such a tree inside a graph of large minimum degree. 

\begin{lemma}[{\cite[Lemma 12]{BDSZZ24b}}]
\label{lm:tree}
Let $G$ be a graph of minimum degree at least $d$. Then $G$ has an induced subgraph $H$ of minimum degree at least $d$ with a spanning tree of maximum degree at most $d$.
\end{lemma}

For a tree $T$ on $m$ vertices with maximum degree at most $\Delta\geq 1$, we have
\begin{equation}
\label{eq:tree-automorphisms}
    |\mathrm{Aut}(T)|\leq m\Delta^{m-1}.
\end{equation}
Indeed, fix a root in $T$ in an arbitrary way. Then, choose the image of the root in at most $m$ ways, then consider the remaining vertices in an order in which each parent precedes its children. Each vertex must map to a neighbour of its parent's image, giving at most $\Delta$ choices. For fixed $\Delta$, the bound is exponential in~$m$.

%%%%%%%%%%%%%%%%%%%%%%%%%%%%%%%%%%
\section{Proof of \texorpdfstring{\cref{thm:small-linear-nc}}{Theorem~\ref{thm:small-linear-nc}}} \label{sec:main}
%%%%%%%%%%%%%%%%%%%%%%%%%%%%%%%%%%

In this section we prove our main result, i.e.,\ we show that every small class has linear neighbourhood complexity. 

Assume towards a contradiction that there exists a small class $\mathcal{C}$ with superlinear neighbourhood complexity. 
Let $c$ be such that $|\cC_n| \leq n!c^n$ holds for all $n \in \bN$. 
Then, for every integer $C$, there exist a graph $G\in\cC$ and a non-empty set $A\subseteq V(G)$ such that 
\[
    |\cN_G(A)|\geq (C+1)|A|.
\]
Note that vertices in $A$ contribute at most $|A|$ neighbourhoods to $\cN_G(A)$. Thus, there exists a set $B\subseteq V(G)\setminus A$ such that $B$ is $A$-diverse and $|B|\geq C|A|$.

Since $\cC$ is hereditary, after passing to an induced subgraph of $G[A \cup B]$ and interchanging the two parts if necessary, we may assume that $|A|+|B|$ is minimal among all pairs with these properties, i.e.,\ they satisfy the following minimality conditions:
\begin{itemize}
    \item[(M1)] for every non-empty proper subset $A'\subset A$, the vertices of $B$ have fewer than $C|A'|$ distinct $A'$-neighbourhoods;
    \item[(M2)] for every non-empty proper subset $B'\subset B$, the vertices of $A$ have fewer than $C|B'|$ distinct $B'$-neighbourhoods.
\end{itemize}

Denote $F:=G[A,B]$. The main idea behind the proof of \cref{thm:small-linear-nc} is to fix $G$ and minimal sets $A$ and $B$ as above and construct, for some integers $a$ and $N$, more than $N!c^N$ distinct labelled bipartite graphs of the form $G[A_0,B_0]$. Here $A_0\subseteq A$ and $B_0\subseteq B$ satisfy $|A_0|=a$ and $|B_0|=N-a$, and we label $A_0$ by $[a]$ and $B_0$ by $\{a+1,\ldots,N\}$. Although these bipartite graphs need not belong to $\cC$, the corresponding induced graphs $G[A_0, B_0]$ do. Since the labels determine the bipartition, distinct labelled bipartite graphs remain distinct after the edges inside $A_0$ and $B_0$ are restored. We thus obtain more than $N!c^N$ labelled graphs in $\cC_N$, contradicting the assumption of smallness of $\cC$. 

\begin{figure}[htbp]
	\centering
	\begin{tikzpicture}
		% Base template for the set-containment diagrams used around Theorem~\ref{thm:auto}.
% This file intentionally contains no tikzpicture environment.  Later figures can
% input it and add annotations using the named nodes
% auto-A, auto-Aprime, auto-B, auto-Bprime, auto-Btilde, and auto-Bdoubleprime.

\node[
	containment/ambient,
	minimum width=.60\linewidth,
	minimum height=16mm
] (auto-A) {};

\node[
	containment/selected,
	draw=black!65,
	fill=FigRed!18,
	minimum width=.40\linewidth,
	minimum height=15mm
] (auto-Aprime) at ($(auto-A.center)+(8mm,0)$) {};

\node[
	containment/ambient,
	minimum width=.82\linewidth,
	minimum height=20mm,
	below=14mm of auto-A
] (auto-B) {};

\node[
	containment/selected,
	draw=black!65,
	fill=FigPoolFill,
	minimum width=.64\linewidth,
	minimum height=19mm
] (auto-Bprime) at ($(auto-B.center)+(2mm,0)$) {};

% An extension of the core chosen in Section~\ref{sec:MainTheoremFromLNC}.
\node[
	containment/selected,
	draw=black!65,
	fill=FigBlue!40,
	line width=.1pt,
	minimum width=.45\linewidth,
	minimum height=15mm
] (auto-Btilde) at ($(auto-Bprime.center)+(10mm,0)$) {};

\node[
	containment/core,
	draw=black!65,
	fill=FigCoreFill,
	minimum width=.30\linewidth,
	minimum height=11mm
] (auto-Bdoubleprime) at ($(auto-Bprime.center)+(18mm,0)$) {};

\node[font=\small, anchor=east]
	at ($(auto-A.west)+(-3mm,0)$) {$A$};
\node[containment/label, font=\fontsize{7.3pt}{8.8pt}\selectfont]
	at ($(auto-Aprime.north west)+(2mm,-2mm)$) {$A'$};
\node[font=\small, anchor=east]
	at ($(auto-B.west)+(-3mm,0)$) {$B$};
\node[containment/label, font=\fontsize{7.3pt}{8.8pt}\selectfont]
	at ($(auto-Bprime.north west)+(2mm,-2mm)$) {$B'$};
\node[containment/label, font=\fontsize{7.3pt}{8.8pt}\selectfont]
	at ($(auto-Btilde.north west)+(2mm,-2mm)$) {$\tilde B$};
\node[containment/label, font=\fontsize{7.3pt}{8.8pt}\selectfont]
	at ($(auto-Bdoubleprime.north west)+(2mm,-2mm)$) {$B''$};
	\end{tikzpicture}
	\caption{The relationship between the sets \setcaptionbox{FigRed!25}{A'}, \setcaptionbox{FigBlue!25}{B'}, and \setcaptionbox{FigOrange!25}{B''} guaranteed by \cref{thm:auto}. The set \setcaptionbox{FigBlue!40}{\tilde B} is an extension of $B''$ inside $B'$.}
	\label{fig:auto-sets}
\end{figure}
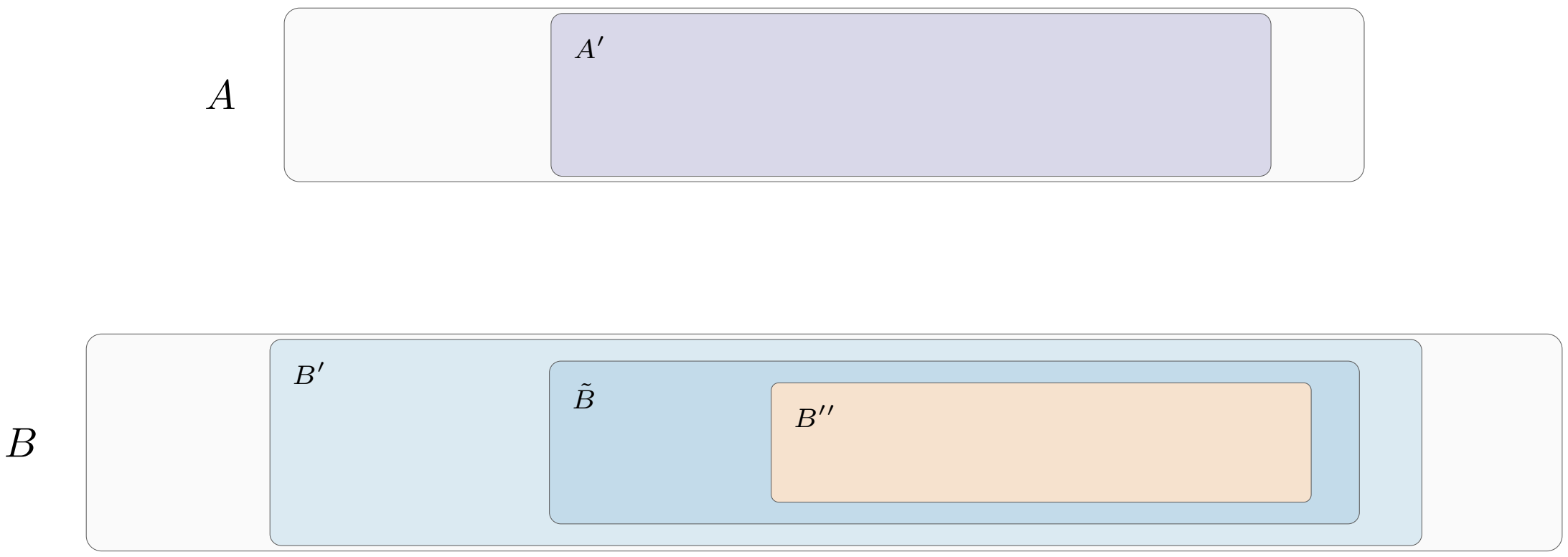

Our main technical ingredient is the following theorem. It finds sets $A'\subseteq A$ and $B''\subseteq B'\subseteq B$ such that $F[A', B'']$ has few automorphisms and $B'\setminus B''$ is $A'$-diverse (see \cref{fig:auto-sets}). The large set $B'$ provides many extensions $\tilde B$ with $B'' \subseteq \tilde B \subset B'$. The $A'$-diversity of $B' \setminus B''$ and the automorphism bound allow us to control how many extensions and labellings can produce the same labelled graph.
\begin{theorem}
\label{thm:auto}
Let $C_0$ be a sufficiently large integer and set $C := C_0^4$.
Let $F$ be a 
bipartite graph with non-empty parts $A,B$ such that $B$ is $A$-diverse, $|B|\geq C|A|$, and the minimality conditions (M1) and (M2) are satisfied. Then, there exist subsets $A'\subseteq A$, $B''\subseteq B'\subseteq B$, and a positive integer $n'$ that satisfy the following properties: 
\begin{itemize}
    \item[(1)]\label{itm:a1} $|A'|\leq n'$, $|B'|=C_0n'$, and $|B''|\leq 3n'$,
    \item[(2)]\label{item:auto2} the set $B'\setminus B''$ is $A'$-diverse,
    \item[(3)]\label{item:auto3} the induced subgraph $F[A',B'']$ has at most $C_0^{100n'}$ automorphisms.
\end{itemize}
\end{theorem}

In the next section, we show how to derive our main result (\cref{thm:small-linear-nc}) from \cref{thm:auto}. Then, in \cref{sec:proof-of-auto}, we prove \cref{thm:auto}.

\subsection[Theorem~\ref{thm:small-linear-nc} follows from Theorem~\ref{thm:auto}]{\cref{thm:small-linear-nc} follows from \cref{thm:auto}}
\label{sec:MainTheoremFromLNC}

Assume \cref{thm:auto}, and choose $C_0>c^3$ to be a sufficiently large perfect square. Let $A',B',B''$ and $n'$ be given by the theorem, and set $a':=|A'|$ and
\[
 N:=a'+\sqrt{C_0}n'.
\]
We will construct at least $N!C_0^{N/3} > N!c^N$ distinct labelled bipartite graphs of the form $F[A',\tilde B]$, with the labels determining the two parts. As explained at the beginning of \cref{sec:main}, restoring the edges within each part then gives at least as many graphs in $\cC_N$, contradicting smallness.

Let $\cF$ be the family of labelled bipartite graphs with two parts $[a']$ and $\{a'+1,\ldots,N\}$ isomorphic to a graph $F[A',\tilde B]$ for some $\tilde B$ with $B''\subseteq\tilde B\subseteq B'$ and $|\tilde B|=\sqrt{C_0}n'$. In order to show that $|\mathcal{F}|> N!C_0^{N/3}$, we will count labellings of graphs $F[A',\tilde B]$, where $A'$ gets labelled by the elements in $[a']$ and $\tilde B$ gets labelled by the elements in $\{a'+1,\ldots,N\}$. 

First, we count the number of ways to select the extension $\tilde B$ of $B''$. Since $|B''|\leq 3n'$ and $|B'|=C_0n'$, there are at least $\binom{(C_0-3)n'}{(\sqrt{C_0}-3)n'}$ such extensions. Once $\tilde B$ is fixed, there are $a'!$ ways to label $A'$ and $(N-a')!$ ways to label $\tilde B$. Note that different choices of $\tilde B$ and the labellings of $A'$ and $\tilde B$ may produce the same labelled graph. Let $K$ be the largest number of choices producing the same labelled graph. Then 
\[
 	|\mathcal{F}|\geq\frac{a'!(N-a')!}{K}\binom{(C_0-3) n'}{(\sqrt{C_0}-3)n'}.
\]
We now bound $K$ from above.
For this, fix a labelled graph $\tilde F$ obtained by labelling $F[A',\tilde B]$ for some $\tilde B$. The choices of extensions and labellings that produce the same labelled graph $\tilde F$ correspond to isomorphisms $f$ from $\tilde F$ onto induced subgraphs of $F[A',B']$ that preserve $A'$ setwise and satisfy $B'' \subseteq f(\tilde B)$. There are at most $\binom{C_0n'}{3n'}$ choices for $f(B'')$. Once this image is fixed, there are at most $|\mathrm{Aut}(F[A',B''])|$ ways to specify $f$ on $A'\cup B''$.

This restriction determines the rest of $f$. Indeed, every vertex of $B''\setminus f(B'')$ must lie in the image, and its preimage is uniquely determined by the distinct neighbourhoods of vertices in $\tilde B\setminus B''$ on $A'$. The remaining vertices must map into $B'\setminus B''$, where distinct neighbourhoods on $A'$ uniquely determine their images. Therefore,
\[
 K\leq\binom{C_0n'}{3n'} C_0^{100n'}.
\]

Combining these bounds and using $N=a'+\sqrt{C_0}n'$ and $a'\leq n'$, we obtain
\begin{align*}
 |\mathcal{F}|&\geq\frac{a'!(N-a')!}{ \binom{C_0n'}{3n'}\cdot C_0^{100n'}}\binom{(C_0-3) n'}{(\sqrt{C_0}-3)n'}\\
 &=\frac{N!}{\binom{N}{a'} \binom{C_0n'}{3n'}\cdot C_0^{100n'}}\binom{(C_0-3) n'}{(\sqrt{C_0}-3)n'}\\
 &\geq\frac{N!}{2^N(eC_0/3)^{3n'}C_0^{100n'}} \sqrt{C_0}^{(\sqrt{C_0}-3)n'}\\
 &\geq N!\left(\frac{C_0^{(\sqrt{C_0}-3)/2}}{2^{\sqrt{C_0}+1}(eC_0/3)^3C_0^{100}}\right)^{n'}
 >N! \cdot C_0^{N/3},
\end{align*}
provided that $C_0$ is sufficiently large. This completes the proof.

\subsection[Proof of Theorem~\ref{thm:auto}]{Proof of \cref{thm:auto}}
\label{sec:proof-of-auto}

By deleting vertices of $B$ if necessary, we may assume $|B|=C|A|$. This preserves $A$-diversity and the minimality conditions (M1) and (M2). Set $n:=|A|$. 

If $n\leq C^{10}$, set $A':=A$, $n':=n$, and $B'':=\emptyset$, and choose any $B'\subseteq B$ of size $C_0n$. Such a choice is possible since $|B|=C_0^4n$. Properties~(1) and~(2) hold immediately. The graph $F[A',B'']$ is edgeless, so
\[
    |\mathrm{Aut}(F[A',B''])|=n!\leq n^n\leq C^{10n}=C_0^{40n}\leq C_0^{100n'},
\]
which gives property~(3). Thus the theorem holds in this case, and we assume $n>C^{10}$ for the rest of the proof.

Recall that \cref{lm:bounded} provides a subgraph $F[A, B']$ with few automorphisms when every vertex of $B$ has a bounded number of neighbours in $A$. This suggests looking for a suitable restriction $F[A',B']$ with $B'$ large and all its $A'$-neighbourhoods bounded in size. We do not know how to obtain such a restriction in general. Instead, we find a large set $B_{\good}\subseteq B$ of vertices whose $A$-neighbourhoods differ from suitable canonical neighbourhoods in a number of vertices bounded in terms of $C$. We call these vertices \emph{good} and the corresponding symmetric differences \emph{defects}.

We use these defects to choose $A',B',B''$ and construct an
auxiliary bipartite graph in which we find a core subgraph with
few automorphisms. Depending on the number of distinct defects,
we obtain this core using either \cref{lm:bounded} or
\cref{lm:tree}. We then construct an injective homomorphism from
a suitable setwise stabiliser in $\mathrm{Aut}(F[A',B''])$ to
the automorphism group of the auxiliary core. This stabiliser
is defined by a fixed number of subsets of $A'\cup B''$, so the orbit--stabiliser theorem implies that the full automorphism group is larger by at most a factor exponential in $|A'|+|B''|$ (\cref{lm:orbit-stabiliser}).
Together with the bound on the number of automorphisms of the
auxiliary core, this gives property~(3).

\subsubsection{Selecting many good vertices from \texorpdfstring{$B$}{B}}
\label{sec:many-good}

We first partition $A$ into equivalence classes and choose one representative from each class. For each vertex $v\in B$, its \emph{canonical neighbourhood} is the union of those classes whose representatives are adjacent to $v$. We call a vertex \emph{good} if its $A$-neighbourhood differs from its canonical neighbourhood in only a bounded number of vertices. The packing lemma and a random choice of representatives will ensure that at least a third of the vertices in $B$ are good.

Let $R\subseteq B$ be an inclusion-maximal subset such that $|N_A(u)\triangle N_A(v)|>C^{10}$ for all distinct $u,v\in R$. Partition $A$ into equivalence classes $A_1,\ldots,A_k$, where two vertices are equivalent if they have the same neighbourhood in~$R$.

\begin{figure}[htbp]
	\centering
	\begin{tikzpicture}
		% Set layout for Lemma~\ref{lm:profile-defect}.
% This file has no tikzpicture environment so later figures can input it
% and add annotations using the named nodes and coordinates below.

\node[
	containment/ambient,
	minimum width=.60\linewidth,
	minimum height=26mm
] (profile-A) {};

\node[
	containment/ambient,
	minimum width=.82\linewidth,
	minimum height=24mm,
	below=14mm of profile-A
] (profile-B) {};

% Show the first two and last classes; the others are omitted.
\foreach \fraction in {0.12,0.24,0.88} {
	\draw[black!65, line width=.1pt]
		($(profile-A.north west)!\fraction!(profile-A.north east)$) --
		($(profile-A.south west)!\fraction!(profile-A.south east)$);
}

% Q contains one representative from every class.  Its unfilled outline
% leaves the partition dividers visible.
\node[
	containment/selected,
	draw=black!65,
	fill=none,
	minimum width=\dimexpr.60\linewidth-2mm\relax,
	minimum height=5mm
] (profile-Q) at ($(profile-A.center)+(0,-6.5mm)$) {};

\node[font=\fontsize{7.3pt}{8.8pt}\selectfont, inner sep=0pt] at
	($(profile-Q.west)!0.40!(profile-Q.east)$) {$Q$};

\node[font=\small, anchor=east] at ($(profile-A.west)+(-3mm,0)$) {$A$};
\node[font=\small, anchor=east] at ($(profile-B.west)+(-3mm,0)$) {$B$};

\node[font=\fontsize{7.3pt}{8.8pt}\selectfont, anchor=north, inner sep=0pt] at
	($(profile-A.north west)!0.06!(profile-A.north east)+(0,-1mm)$) {$A_1$};
\node[font=\fontsize{7.3pt}{8.8pt}\selectfont, anchor=north, inner sep=0pt] at
	($(profile-A.north west)!0.18!(profile-A.north east)+(0,-1mm)$) {$A_2$};
\node[font=\fontsize{7.3pt}{8.8pt}\selectfont, anchor=north, inner sep=0pt] at
	($(profile-A.north west)!0.94!(profile-A.north east)+(0,-1mm)$) {$A_k$};
\node[font=\fontsize{7.3pt}{8.8pt}\selectfont, inner sep=0pt] at
	($(profile-A.north west)!0.56!(profile-A.north east)+(0,-2.3mm)$)
	{$\cdots$};

\node[circle, fill=black!70, inner sep=0pt, minimum size=1.2mm]
	(profile-q-one) at
	($(profile-A.south west)!0.06!(profile-A.south east)+(0,6.5mm)$) {};
\node[font=\fontsize{7.3pt}{8.8pt}\selectfont, anchor=west, inner sep=0pt] at
	($(profile-q-one.east)+(1mm,0)$) {$q_1$};
\node[circle, fill=black!70, inner sep=0pt, minimum size=1.2mm]
	(profile-q-two) at
	($(profile-A.south west)!0.18!(profile-A.south east)+(0,6.5mm)$) {};
\node[font=\fontsize{7.3pt}{8.8pt}\selectfont, anchor=west, inner sep=0pt] at
	($(profile-q-two.east)+(1mm,0)$) {$q_2$};
\node[circle, fill=black!70, inner sep=0pt, minimum size=1.2mm]
	(profile-q-last) at
	($(profile-A.south west)!0.94!(profile-A.south east)+(0,6.5mm)$) {};
\node[font=\fontsize{7.3pt}{8.8pt}\selectfont, anchor=west, inner sep=0pt] at
	($(profile-q-last.east)+(1mm,0)$) {$q_k$};
\node[font=\fontsize{7.3pt}{8.8pt}\selectfont, inner sep=0pt] at
	($(profile-A.south west)!0.56!(profile-A.south east)+(0,6.5mm)$)
	{$\cdots$};

\node[
	containment/ambient,
	minimum width=.19\linewidth,
	minimum height=23mm
] (profile-R) at ($(profile-B.center)+(-.215\linewidth,0)$) {};

\node[
	containment/ambient,
	minimum width=.43\linewidth,
	minimum height=23mm
] (profile-Bgood) at ($(profile-B.center)+(.105\linewidth,0)$) {};

\node[containment/label, font=\fontsize{7.3pt}{8.8pt}\selectfont] at
	($(profile-Bgood.north west)+(2mm,-2mm)$) {$B_{\good}$};
\node[containment/label, font=\fontsize{7.3pt}{8.8pt}\selectfont] at
	($(profile-R.north west)+(2mm,-2mm)$) {$R$};
	\end{tikzpicture}
	\caption{Illustration for \cref{lm:profile-defect}: the partition $A=A_1\cup\cdots\cup A_k$, the representative set $Q=\{q_1,\ldots,q_k\}$ with $q_i\in A_i$, and the disjoint sets $B_{\good}$ and $R$ in $B$.}
	\label{fig:profile-defect-sets}
\end{figure}
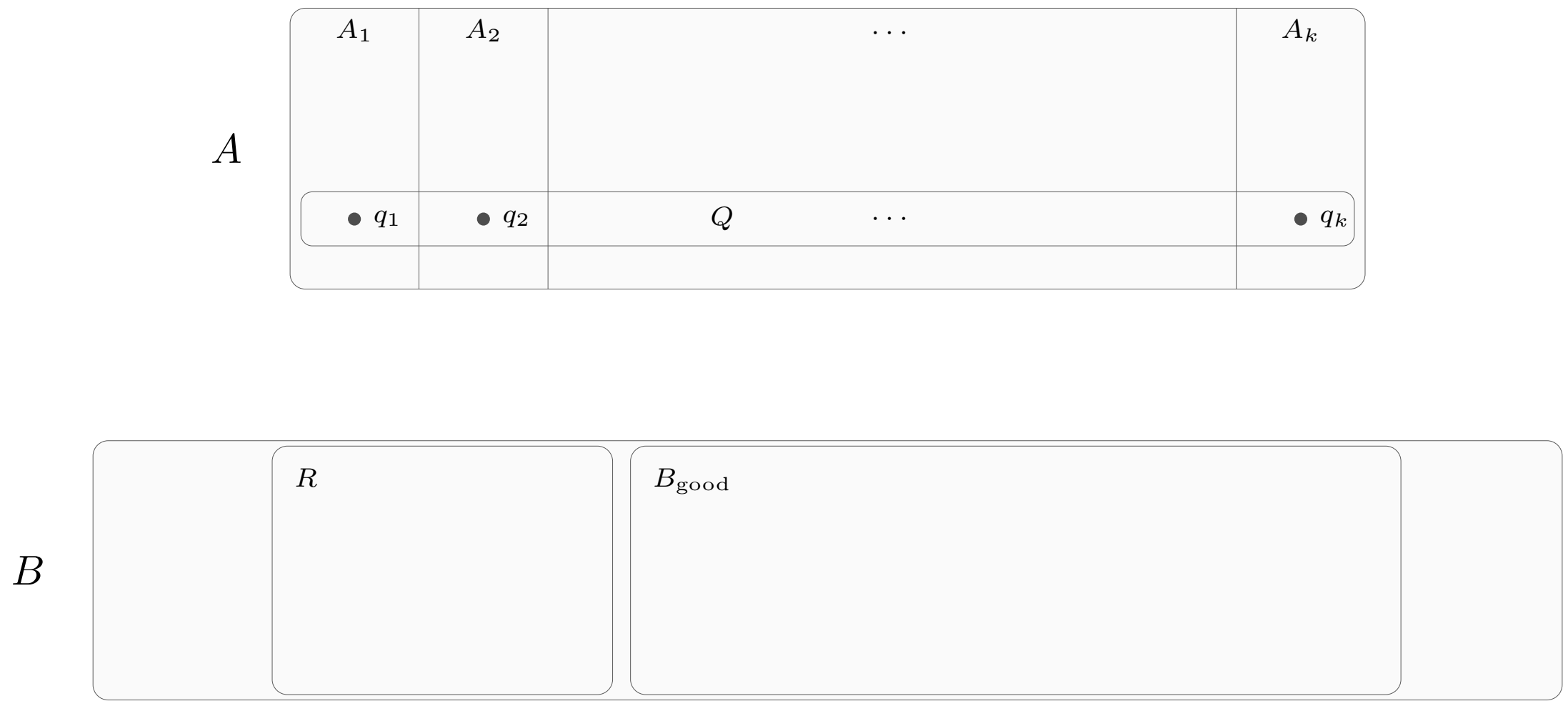

\begin{lemma}
\label{lm:profile-defect}
There exist a set $B_{\good}\subseteq B\setminus R$ of size at least $|B|/3$ and a set $Q=\{q_1,\ldots,q_k\}\subseteq A$ of representatives $q_i\in A_i$ such that the following holds. For every vertex $v\in B$, let $\Sigma(v)\subseteq A$ be the union of all equivalence classes $A_i$ such that $q_i$ is adjacent to $v$ in $F$. Then, for every $v\in B_{\good}$, $|N_A(v)\triangle\Sigma(v)|\leq 4C^{10}$.
\end{lemma}
\begin{proof}
See \cref{fig:profile-defect-sets} for an illustration. 
Recall that $n>C^{10}$.
By \cref{cor:packing}, $|R|\leq n/C^7$. For every $v\in B\setminus R$, find a vertex $r(v)\in R$ such that $|N_A(v)\triangle N_A(r(v))|\leq C^{10}$ --- this is possible due to the maximality of $R$. For all $i\in[k]$, take $q_i\in A_i$ uniformly at random. Fix $v\in B\setminus R$. For $i\in[k]$, let $D_i(v)\subseteq A_i$ be the set of vertices that have different adjacencies with $v$ and $r(v)$. We know that $\sum_{i=1}^k|D_i(v)|\leq C^{10}$. Moreover, all vertices in $A_i$ have the same adjacency relation with $r(v)$. Therefore, the following two cases are possible. If $q_i\notin D_i(v)$, then $|N_{A_i}(v)\triangle(\Sigma(v)\cap A_i)|=|D_i(v)|$. Otherwise, $|N_{A_i}(v)\triangle(\Sigma(v)\cap A_i)|\leq|A_i|$. Then
$$
\mathbb{E}|N_A(v)\triangle \Sigma(v)|\leq\sum_{i=1}^k\left(\left(1-\frac{|D_i(v)|}{|A_i|}\right)|D_i(v)|+\frac{|D_i(v)|}{|A_i|}\cdot|A_i|\right)\leq 2\sum_{i=1}^k|D_i(v)|\leq 2C^{10}.
$$
Therefore, by Markov's inequality,
$$
 \mathbb{P}(|N_A(v)\triangle \Sigma(v)|>4C^{10})\leq 1/2.
$$
Let $\xi$ be the number of vertices $v\in B\setminus R$ that have $|N_A(v)\triangle \Sigma(v)|>4C^{10}$. We get $\mathbb{E}\xi\leq |B\setminus R|/2$. Therefore, $\mathbb{P}(\xi\leq|B\setminus R|/2)>0$. So, there exist $q_1,\ldots,q_k$ such that at least $|B\setminus R|/2$ vertices $v$ have $|N_A(v)\triangle \Sigma(v)|\leq 4C^{10}$. We get $|B\setminus R|/2\geq |B|/2-|R|>|B|/3$, completing the proof.
\end{proof}

Fix sets $B_{\good}$ and $Q$ as in \cref{lm:profile-defect}, so $|B_{\good}|\geq |B|/3$. For each $v\in B_{\good}$, define
\begin{itemize}
\item its \emph{profile} $P_v:=N_Q(v)$, and
\item its \emph{defect} $D_v:=N_A(v)\triangle\Sigma(v)$.
\end{itemize}
Every defect has size at most $4C^{10}$ and is disjoint from $Q$, since $N_A(v)$ and $\Sigma(v)$ agree on the representatives. Let $d:=|\{D_v:\,v\in B_{\good}\}|$ be the number of distinct defects. The profile $P_v$ determines $\Sigma(v)$, and $N_A(v)=\Sigma(v)\triangle D_v$. Since $B_{\good}$ is $A$-diverse, we obtain the following observation.

\begin{observation}
\label{obs:profile-defect}
	For any distinct vertices $u,v \in B_{\good}$, we have $(P_u,D_u) \neq (P_v,D_v)$.
\end{observation}

For $u\in A_i$, define its \emph{representative} by $\rep(u):=q_i$. For $A_0\subseteq A$, write $\rep(A_0):=\{\rep(u):\,u\in A_0\}$. We call a set $R_0\subseteq R$ a \emph{witness} for $A_0$ if $\rep(A_0)$ is $R_0$-diverse. Thus, the neighbourhood of a vertex of $A_0$ in $R_0$ uniquely identifies its representative within $\rep(A_0)$. The next lemma shows that an automorphism preserving the four sets in its statement setwise sends each vertex's representative to the representative of its image. It also preserves defect membership within the selected sets, since membership is determined by whether a good vertex has different adjacencies to a vertex and its representative.

\begin{lemma}
\label{lem:rep}
Let $A_0\subseteq A\setminus Q$ and $B_0\subseteq B_{\good}$. Let $Q_0=\rep(A_0)$. Let $R_0\subseteq R$ be a witness for $A_0$. Let $\sigma$ be an automorphism of the graph $F[A_0\cup Q_0, B_0\cup R_0]$ that preserves the sets $A_0,B_0,Q_0,R_0$. Then, for every $u\in A_0$, we have $\sigma(\rep(u))=\rep(\sigma(u))$. Moreover, for every $u\in A_0$ and $v \in B_0$, we have $u\in D_v$ if and only if $\sigma(u)\in D_{\sigma(v)}$.
\end{lemma}

\begin{proof}
Fix $\sigma$ and $u\in A_0$, as in the statement of the lemma. We have that $N_R(u)=N_R(\rep(u))$ and therefore $N_{R_0}(u)=N_{R_0}(\rep(u))$. Since $\sigma$ is an automorphism of $F[A_0\cup Q_0, B_0\cup R_0]$ that preserves $R_0$, we get that $N_{R_0}(\sigma(u))=N_{R_0}(\sigma(\rep(u)))$. On the other hand, by the definition of representatives, $N_R(\sigma(u))=N_R(\rep(\sigma(u)))$ and therefore $N_{R_0}(\sigma(u))=N_{R_0}(\rep(\sigma(u)))$. Since all vertices in $Q_0$ have different neighbourhoods in $R_0$ and since $\sigma$ preserves $Q_0$, we get $\rep(\sigma(u)),\sigma(\rep(u))\in Q_0$, implying $\rep(\sigma(u))=\sigma(\rep(u))$.

For the second part of the lemma, observe
\begin{align*}
 u\in D_v&\Leftrightarrow v\text{ has different adjacency relations with }u\text{ and }\rep(u)\\
 &\Leftrightarrow \sigma(v)\text{ has different adjacency relations with }\sigma(u)\text{ and }\sigma(\rep(u))\\
  &\Leftrightarrow \sigma(v)\text{ has different adjacency relations with }\sigma(u)\text{ and }\rep(\sigma(u))\\
  &\Leftrightarrow \sigma(u)\in D_{\sigma(v)}.
\end{align*}
This completes the proof. 
\end{proof}

We now construct $A'$, $B'$, and $B''$ in two cases, according to the number $d$ of distinct defects. If $d\geq\sqrt{C}n$ (\cref{sc:many}), we use the auxiliary bipartite graph $F_{\mathrm{defect}}$ with parts $A$ and $B_{\good}$ in which $u$ and $v$ are adjacent if and only if $u\in D_v$. Selecting one vertex for each distinct defect gives a large diverse set with bounded neighbourhoods, allowing us to apply \cref{lm:bounded} after passing to a minimal restriction. This yields a core subgraph of $F_{\mathrm{defect}}$ with few automorphisms.

If $d<\sqrt{C}n$ (\cref{sc:few}), we instead use the auxiliary bipartite graph $F_{\prodef}$ whose parts are disjoint copies of the families of profiles and defects, with an edge between $P$ and $D$ whenever $(P,D)=(P_v,D_v)$ for some $v\in B_{\good}$. By \cref{obs:profile-defect}, this graph has $|B_{\good}|\geq Cn/3$ edges. The packing bound and the minimality conditions also bound the number of profiles, so the graph has few vertices relative to its number of edges. This allows us to find a bounded-degree tree $T$ with few automorphisms, which serves as a core subgraph of $F_{\prodef}$.

In both cases, we include suitable representative and witness sets $Q_0,R_0$ in $A'$ and $B''$, respectively. We then use \cref{lem:rep} to show that a suitable setwise stabiliser in $\mathrm{Aut}(F[A',B''])$ preserves defect incidences, allowing us to construct an injective homomorphism into the automorphism group of the core chosen in the corresponding auxiliary graph.

\begin{remark}\label{rem:B'-truncation}
    In both cases, it suffices to prove property~(1) with $|B'|\geq C_0n'$ in place of $|B'|=C_0n'$. Once the other requirements have been verified, we can reduce $B'$ to size $C_0n'$ by deleting vertices from $B'\setminus B''$, since $|B''|\leq 3n'\leq C_0n'$. This preserves the $A'$-diversity required in property~(2) and leaves the graph $F[A',B'']$ in property~(3) unchanged.
\end{remark}

\subsubsection{Many different defects: \texorpdfstring{$d\geq\sqrt{C}n$}{d at least square root C times n}}
\label{sc:many}

\paragraph{Auxiliary graph.}
Define the auxiliary bipartite graph $F_{\mathrm{defect}}$ with parts $A$ and $B_{\good}$ where every vertex $v\in B_{\good}$ is adjacent to every vertex in $D_v$. Since $d\geq\sqrt{C}n$, by \cref{lm:minimal-diverse}, there exist a non-empty set $A_0\subset A$ and a minimally $A_0$-diverse set $B_0\subset B_{\good}$ in $F_{\mathrm{defect}}$ such that $|B_0|=\sqrt{C}|A_0|$.

Recall that $D_v\subset A\setminus Q$ for all $v\in B_{\good}$. Therefore, $A_0\subset A\setminus Q$. Indeed, otherwise the graph $F_{\mathrm{defect}}[A_0,B_0]$ would have an isolated vertex in $A_0$, which would contradict the minimality condition.

\paragraph{Core of the auxiliary graph.}
Since every vertex in $B_{\good}$ has at most $4C^{10}$ neighbours in $F_{\mathrm{defect}}$, applying \cref{lm:bounded} to the induced subgraph $F_{\mathrm{defect}}[A_0,B_0]$, we conclude that there exists a set $B'_0\subset B_0$ of size at most $|A_0|-1$ such that the group of automorphisms $\mathrm{Aut}(F_{\mathrm{defect}}[A_0,B'_0])$ has size at most exponential in $|A_0|$. The graph $F_{\mathrm{defect}}[A_0,B'_0]$ is the core of $F_{\mathrm{defect}}$ that we will use to establish property~(3).

\paragraph{The choice of \texorpdfstring{$A',B',B''$ and $n'$}{A',B',B'' and  n'}.} Let us get back to the original graph $F$. Let $Q_0=\rep(A_0)$. Recall that vertices in $Q_0\subset Q$ have different neighbourhoods in $R$. Due to \cref{lm:small-diverse}, we can choose $R_0\subset R$ such that 
\begin{itemize}
\item $|R_0|\leq|Q_0|-1$,
\item vertices in $R_0$ have different traces on $Q_0$ in $F$,
\item $R_0$ is a witness of $A_0$, that is vertices in $Q_0$ have different traces on $R_0$ in $F$.
\end{itemize}
In particular, for every $i\in[k]$ and every $x\in A_i\cap A_0$ the vertex $q_i$ is the unique vertex in $Q_0$ that has exactly the same adjacency relations with $R_0$ in $F$.

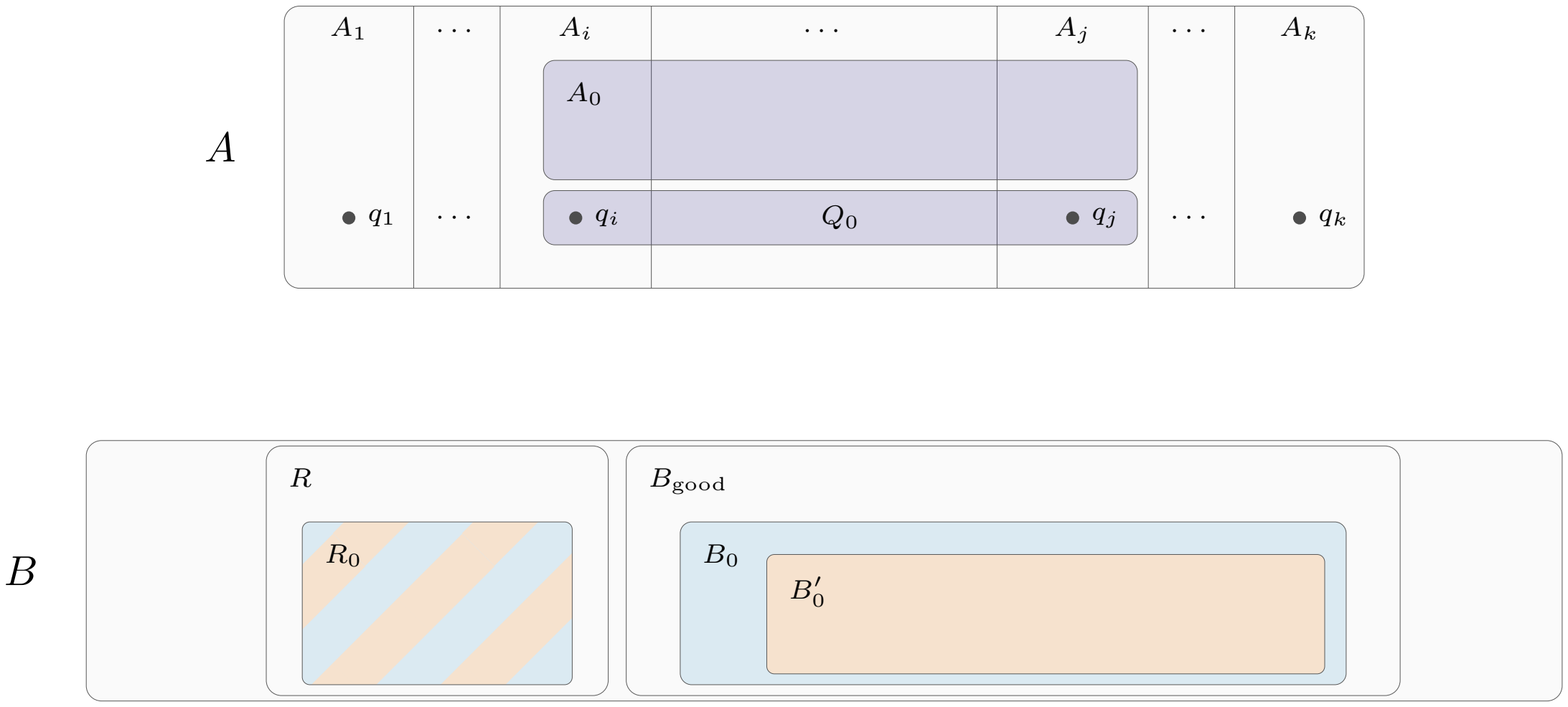
\begin{figure}[htbp]
	\centering
	\begin{tikzpicture}
		% Set layout for the many-defects case of Theorem~\ref{thm:auto}.
% This file has no tikzpicture environment, so its named nodes can be
% reused and annotated in later figures.

\node[
	containment/ambient,
	minimum width=.60\linewidth,
	minimum height=26mm
] (many-A) {};

% The displayed classes are A_1, A_i, A_j, and A_k, with omitted classes
% indicated by ellipses.  The dividers remain visible through A_0 and Q_0.
\foreach \fraction in {0.12,0.20,0.34,0.66,0.80,0.88} {
	\draw[black!65, line width=.1pt]
		($(many-A.north west)!\fraction!(many-A.north east)$) --
		($(many-A.south west)!\fraction!(many-A.south east)$);
}

% A_0 begins inside A_i and ends inside A_j; Q_0 lies below it.
% Their translucent fills leave the partition dividers visible.
\node[
	containment/selected,
	draw=black!65,
	fill=FigRed,
	fill opacity=.18,
	minimum width=.33\linewidth,
	minimum height=11mm
] (many-Azero) at ($(many-A.center)+(.009\linewidth,2.5mm)$) {};

\node[
	containment/selected,
	draw=black!65,
	fill=FigRed,
	fill opacity=.18,
	minimum width=.33\linewidth,
	minimum height=5mm
] (many-Qzero) at ($(many-A.center)+(.009\linewidth,-6.5mm)$) {};

\node[
	containment/ambient,
	minimum width=.82\linewidth,
	minimum height=24mm,
	below=14mm of many-A
] (many-B) {};

% R and B_good are disjoint.  R_0 occurs in both B' and B''.
\node[
	containment/ambient,
	minimum width=.19\linewidth,
	minimum height=23mm
] (many-R) at ($(many-B.center)+(-.215\linewidth,0)$) {};

\node[
	containment/ambient,
	minimum width=.43\linewidth,
	minimum height=23mm
] (many-Bgood) at ($(many-B.center)+(.105\linewidth,0)$) {};

\node[
	containment/core,
	draw=none,
	outer sep=0pt,
	fill=FigPoolFill,
	minimum width=.15\linewidth,
	minimum height=15mm
] (many-Rzero) at ($(many-R.center)+(0,-3mm)$) {};

% Equal blue and peach bands: R_0 belongs to both unions.
% Draw the outline after the striped fill to keep the border continuous.
\path[containment/overlap, draw=black!65, line width=.1pt]
	(many-Rzero.north west) rectangle (many-Rzero.south east);

\node[
	containment/selected,
	draw=black!65,
	fill=FigPoolFill,
	minimum width=.37\linewidth,
	minimum height=15mm
] (many-Bzero) at ($(many-Bgood.center)+(0,-3mm)$) {};

\node[
	containment/core,
	draw=black!65,
	fill=FigCoreFill,
	minimum width=.31\linewidth,
	minimum height=11mm
] (many-Bzeroprime) at ($(many-Bzero.center)+(3mm,-1mm)$) {};

\node[font=\small, anchor=east] at ($(many-A.west)+(-3mm,0)$) {$A$};
\node[font=\small, anchor=east] at ($(many-B.west)+(-3mm,0)$) {$B$};

\node[font=\fontsize{7.3pt}{8.8pt}\selectfont, anchor=north, inner sep=0pt] at
	($(many-A.north west)!0.06!(many-A.north east)+(0,-1mm)$) {$A_1$};
\node[font=\fontsize{7.3pt}{8.8pt}\selectfont, anchor=north, inner sep=0pt] at
	($(many-A.north west)!0.27!(many-A.north east)+(0,-1mm)$) {$A_i$};
\node[font=\fontsize{7.3pt}{8.8pt}\selectfont, anchor=north, inner sep=0pt] at
	($(many-A.north west)!0.73!(many-A.north east)+(0,-1mm)$) {$A_j$};
\node[font=\fontsize{7.3pt}{8.8pt}\selectfont, anchor=north, inner sep=0pt] at
	($(many-A.north west)!0.94!(many-A.north east)+(0,-1mm)$) {$A_k$};
\foreach \fraction in {0.16,0.50,0.84} {
	\node[font=\fontsize{7.3pt}{8.8pt}\selectfont, inner sep=0pt] at
		($(many-A.north west)!\fraction!(many-A.north east)+(0,-2.3mm)$)
		{$\cdots$};
}

\node[containment/label, font=\fontsize{7.3pt}{8.8pt}\selectfont] at
	($(many-Azero.north west)+(2mm,-2mm)$) {$A_0$};
\node[font=\fontsize{7.3pt}{8.8pt}\selectfont, inner sep=0pt] at
	(many-Qzero.center) {$Q_0$};

% Shown representatives: q_i,q_j belong to Q_0; q_1,q_k do not.
\node[circle, fill=black!70, inner sep=0pt, minimum size=1.2mm]
	(many-q-one) at
	($(many-A.south west)!0.06!(many-A.south east)+(0,6.5mm)$) {};
\node[font=\fontsize{7.3pt}{8.8pt}\selectfont, anchor=west, inner sep=0pt] at
	($(many-q-one.east)+(1mm,0)$) {$q_1$};
\node[circle, fill=black!70, inner sep=0pt, minimum size=1.2mm]
	(many-q-i) at
	($(many-A.south west)!0.27!(many-A.south east)+(0,6.5mm)$) {};
\node[font=\fontsize{7.3pt}{8.8pt}\selectfont, anchor=west, inner sep=0pt] at
	($(many-q-i.east)+(1mm,0)$) {$q_i$};
\node[circle, fill=black!70, inner sep=0pt, minimum size=1.2mm]
	(many-q-j) at
	($(many-A.south west)!0.73!(many-A.south east)+(0,6.5mm)$) {};
\node[font=\fontsize{7.3pt}{8.8pt}\selectfont, anchor=west, inner sep=0pt] at
	($(many-q-j.east)+(1mm,0)$) {$q_j$};
\node[circle, fill=black!70, inner sep=0pt, minimum size=1.2mm]
	(many-q-last) at
	($(many-A.south west)!0.94!(many-A.south east)+(0,6.5mm)$) {};
\node[font=\fontsize{7.3pt}{8.8pt}\selectfont, anchor=west, inner sep=0pt] at
	($(many-q-last.east)+(1mm,0)$) {$q_k$};
\foreach \fraction in {0.16,0.84} {
	\node[font=\fontsize{7.3pt}{8.8pt}\selectfont, inner sep=0pt] at
		($(many-A.south west)!\fraction!(many-A.south east)+(0,6.5mm)$)
		{$\cdots$};
}

\node[containment/label, font=\fontsize{7.3pt}{8.8pt}\selectfont] at
	($(many-R.north west)+(2mm,-2mm)$) {$R$};
\node[containment/label, font=\fontsize{7.3pt}{8.8pt}\selectfont] at
	($(many-Rzero.north west)+(2mm,-2mm)$) {$R_0$};
\node[containment/label, font=\fontsize{7.3pt}{8.8pt}\selectfont] at
	($(many-Bgood.north west)+(2mm,-2mm)$) {$B_{\good}$};
\node[containment/label, font=\fontsize{7.3pt}{8.8pt}\selectfont] at
	($(many-Bzero.north west)+(2mm,-2mm)$) {$B_0$};
\node[containment/label, font=\fontsize{7.3pt}{8.8pt}\selectfont] at
	($(many-Bzeroprime.north west)+(2mm,-2mm)$) {$B'_0$};
	\end{tikzpicture}
	\caption{The choice of \setcaptionbox{FigRed!25}{A'=A_0\cup Q_0}, \setcaptionbox{FigBlue!25}{B'=B_0\cup R_0}, and \setcaptionbox{FigOrange!25}{B''=B'_0\cup R_0} in the many-defects case.}
	\label{fig:many-defects-sets}
\end{figure}

We are now ready to define the sets $A',B',B''$ and the number $n'$ as follows (see \cref{fig:many-defects-sets}):
\[
     A'=A_0\cup Q_0, \quad B'=B_0\cup R_0,\quad B''=B_0'\cup R_0, \quad n'=2|A_0|.
\]
In the rest of this section we verify that these sets and the number satisfy the three properties from \cref{thm:auto}.

\paragraph{Property (1).}
Using the assumption that $C_0$ is large enough, we obtain
$$
    |A'| = |A_0| + |Q_0| \leq 2|A_0| = n',
$$ 
$$
    |B'|\geq|B_0|\geq\sqrt{C}|A_0|\geq\frac{\sqrt{C}}{2}n' = \frac{C_0^2}{2}n' \geq C_0n',
$$
and
$$
    |B''| = |B'_0|+|R_0|\leq|A_0|+|Q_0|=|A'| \leq 3n'.
$$ 
Recalling \cref{rem:B'-truncation}, these establish Property (1).

\paragraph{Property (2).} In order to establish that $B'\setminus B''$ is $A'$-diverse in $F$, we show that its superset $B_0$ is $A'$-diverse in $F$. 
Indeed, let $u,v\in B_0$ be distinct. Since $B_0$ is $A_0$-diverse in
$F_{\mathrm{defect}}$, we have $D_u\cap A_0\neq D_v\cap A_0$.
Suppose $N_{A'}(u)=N_{A'}(v)$, that is $N_{A_0}(u) = N_{A_0}(v)$ and $N_{Q_0}(u)=N_{Q_0}(v)$. The latter implies $\Sigma(u)\cap A_0=\Sigma(v)\cap A_0$. This together with the fact that $D_u\cap A_0=N_{A_0}(u)\triangle(\Sigma(u)\cap A_0)$ and $D_v\cap A_0=N_{A_0}(v)\triangle(\Sigma(v)\cap A_0)$, implies $N_{A_0}(u)\neq N_{A_0}(v)$ --- a contradiction.

\paragraph{Property (3).}
Let $H:=\mathrm{Aut}(F[A', B''])$ and $H_{\mathrm{defect}}:=\mathrm{Aut}(F_{\mathrm{defect}}[A_0, B'_0])$ act on $A'\cup B''$ and $A_0\cup B'_0$, respectively. 
Consider stabilisers
\[
    H_0:=H_{A_0,Q_0,R_0,B_0'} \qquad \qquad (H_{\defect})_0:=(H_{\defect})_{A_0,B_0'}.
\]

Then by \cref{lm:bounded}, we have
\[
    |(H_{\mathrm{defect}})_0|\leq|H_{\mathrm{defect}}|\leq \max\{\sqrt{C}, 4C^{10}\}^{4|A_0|} \leq   C^{21n'}.
\] 
By \cref{lm:orbit-stabiliser}, we have $|H|\leq 2^{4(|A'| + |B''|)}\cdot |H_0|\leq 16^{4n'}\cdot |H_0|$.
Therefore, in order to get property (3) of \cref{thm:auto}, it suffices to show that there exists an injective homomorphism $\varphi:\, H_0\hookrightarrow(H_{\mathrm{defect}})_0$. Indeed, this would imply the desired bound
\[
    |H| 
    \leq 16^{4n'} \cdot |H_0|
    \leq 16^{4n'} \cdot |(H_{\defect})_{0}|
    \leq 16^{4n'} \cdot C^{21n'}
    \leq C^{25n'}
    = C_0^{100n'}.
\]

Since automorphisms in $H_0$ preserve $A_0$ and $B'_0$, we can define the homomorphism $\varphi$ in the natural way as the projection on $A_0\cup B'_0$: for $\sigma\in H_0$, we let $\varphi(\sigma)(v)=\sigma(v)$ for $v\in A_0\cup B'_0$.
Then $\varphi$ is a homomorphism to a group of permutations acting on $A_0\cup B'_0$. Also, since $\sigma$ preserves $A_0$ and $B_0'$ setwise, the same holds true for $\varphi(\sigma)$. In order to show that $\varphi$ is a homomorphism between the groups of automorphisms, we need to prove that $\varphi(\sigma)\in (H_{\defect})_0$ for every $\sigma\in H_0$. Also, $\varphi$ is injective if and only if it has a trivial kernel. In order to establish the latter,
it suffices to show that only the trivial permutation in $H_0$ fixes all vertices in $A_0\cup B'_0$, since for any permutation $\sigma$ that acts non-trivially on $A_0\cup B'_0$ we have that its projection $\varphi(\sigma)$ is non-trivial as well.\footnote{Actually an even stronger result holds: $(H_{Q_0,B'_0})_{(A_0)}$ is trivial.} Therefore, the following claim completes the proof of Property (3).

\begin{claim}
\label{cl:defect-restriction}
$\varphi$ is a homomorphism and $(H_{Q_0})_{(A_0\cup B'_0)}$ is trivial.
\end{claim}
\begin{proof}
Let $\sigma\in H_0$, $v\in B'_0$, and $u\in A_0$. Assume $\{u,v\}\in E(F_{\mathrm{defect}})$. We need to show that $\{\varphi(\sigma)(u),\varphi(\sigma)(v)\}=\{\sigma(u),\sigma(v)\}\in E(F_{\mathrm{defect}})$. 
 According to the definition of $F_{\mathrm{defect}}$, we have that $v$ has different adjacency relations with $\rep(u)$ and $u$. Therefore, $\sigma(v)$ has different adjacency relations with $\sigma(\rep(u))$ and $\sigma(u)$. Due to \cref{lem:rep}, we have $\sigma(\rep(u))=\rep(\sigma(u))$. Thus, $\{\sigma(u),\sigma(v)\}\in E(F_{\mathrm{defect}})$, as required.

 It remains to show that $(H_{Q_0})_{(A_0\cup B'_0)}$ is trivial. Assume $\sigma\in (H_{Q_0})_{(A_0\cup B'_0)}$. Due to \cref{lem:rep}, $\sigma$ also fixes $Q_0$ pointwise. Indeed, for any $q\in Q_0$, there exists $u\in A_0$ such that $q=\rep(u)$. Then, $\sigma(q)=\sigma(\rep(u))=\rep(\sigma(u))=\rep(u)=q$. 
  Finally, since $\sigma$ fixes $R_0$ setwise and since $R_0$ is $Q_0$-diverse, $\sigma$ fixes $R_0$ pointwise. Therefore, $\sigma$ is the identity. This completes the proof.
\end{proof}

\subsubsection{Few different defects: \texorpdfstring{$d<\sqrt{C}n$}{d less than square root C times n}}
\label{sc:few}

\paragraph{Auxiliary graph.}
Recall that the profile and defect of a good vertex $v\in B_{\good}$ are $P_v=N_Q(v)$ and $D_v=N_A(v)\triangle\Sigma(v)$, respectively. We encode the profile--defect pairs realised by good vertices in the auxiliary bipartite graph $F_{\prodef}$. Its parts are disjoint copies of the families $\mathcal{P}:=\{P_v:\,v\in B_{\good}\}$ and $\mathcal{D}:=\{D_v:\,v\in B_{\good}\}$, with an edge $\{P_v,D_v\}$ for each $v\in B_{\good}$.
By \cref{obs:profile-defect}, different good vertices give different edges. Thus $F_{\prodef}$ is simple, and its edges are in bijection with $B_{\good}$, giving
\[
    |E(F_{\prodef})|=|B_{\good}|\geq |B|/3=Cn/3.
\]

To bound the number of vertices of $F_{\prodef}$, recall that $Q$ is $R$-diverse and $|R|\leq n/C^7$ by \cref{cor:packing}. The minimality condition (M2) therefore gives $|Q|<C|R|\leq n/C^6$. In particular, $Q$ is a proper subset of $A$, so (M1) bounds the number of profiles by $|\mathcal{P}|<C|Q|<n/C^5$. Since $|\mathcal{D}|=d<\sqrt{C}n$, we obtain
\[
    |V(F_{\prodef})|=|\mathcal{D}|+|\mathcal{P}|<d+n/C^5<2\sqrt{C}n.
\]

\paragraph{Core of the auxiliary graph.}
First, iteratively delete from $F_{\prodef}$ vertices whose current degree is at most $C^{1/3}$. This removes at most $C^{1/3}|V(F_{\prodef})|<2C^{5/6}n$ edges. For sufficiently large $C$, the remaining graph therefore has at least $Cn/3-2C^{5/6}n\geq Cn/4$ edges and minimum degree at least $\Delta:=\lceil C^{1/3}\rceil$. By \cref{lm:tree}, it contains an induced subgraph $F_{\prodef}[\mathcal{P}_0,\mathcal{D}_0]$, where $\mathcal{P}_0\subseteq\mathcal{P}$ and $\mathcal{D}_0\subseteq\mathcal{D}$, with minimum degree at least $\Delta$ and a spanning tree $T$ of maximum degree at most $\Delta$. By \eqref{eq:tree-automorphisms},
\[
    |\mathrm{Aut}(T)|\leq |V(T)|\Delta^{|V(T)|-1}<C^{|V(T)|},
\]
where the last inequality holds for sufficiently large $C$. The tree $T$ is the core of $F_{\prodef}$ that we will use to establish property~(3).

\paragraph{The choice of \texorpdfstring{$A',B',B''$ and $n'$}{A',B',B'' and  n'}.} Let $B_0$ be the set of vertices $v\in B_{\good}$ that have $P_v\in \mathcal{P}_0$ and $D_v\in\mathcal{D}_0$, that is $\{P_v,D_v\}\in E(F_{\prodef}[\mathcal{P}_0,\mathcal{D}_0])$. 
Thus, we have 
$$
    |B_0| \geq |E(F_{\prodef}[\mathcal{P}_0, \mathcal{D}_0])| \geq C^{1/3}(|\cP_0| + |\cD_0|)/2.
$$
For every $P\in\mathcal{P}_0$, choose an arbitrary $v=v(P)\in B_0$ such that $P = P_v$. In particular, $D_v\in\mathcal{D}_0$. Similarly, for every $D\in\mathcal{D}_0$, choose an arbitrary $u=u(D)\in B_0$ such that $D=D_u$ (and therefore $P_u\in\mathcal{P}_0$). Set 
$$
B_0^{\pro}:=\{v(P):\,P\in\mathcal{P}_0\}\quad\text{ and }\quad B_0^{\ddef}:=\{u(D):\,D\in\mathcal{D}_0\}.
$$
Finally, let $B_0^{\tree}$ be the set of all vertices $v\in B_{\good}$ that have $\{P_v,D_v\}\in E(T)$. 

Now, we choose two subsets of $A$ whose union will be used to distinguish vertices in $B_0$. 
 Due to \cref{lm:small-diverse}, there exist sets $A_0^{\pro}\subseteq Q$ and $A_0^{\ddef}\subseteq A\setminus Q$ such that 
 \begin{itemize}
 \item $|A_0^{\pro}|\leq|\mathcal{P}_0|-1$, $|A_0^{\ddef}|\leq|\mathcal{D}_0|-1$, 
 \item all sets $P\cap A_0^{\pro}$, $P\in\mathcal{P}_0$, and $D\cap A_0^{\ddef}$, $D\in\mathcal{D}_0$, are different, and
 \item all sets $\{P\in\mathcal{P}_0:\,v\in P\}$, $v\in A_0^{\pro}$, are different, and all sets $\{D\in\mathcal{D}_0:\,u\in D\}$, $u\in A_0^{\ddef}$, are different, as well. 
 \end{itemize} 

Let $Q_0:=\rep(A_0^{\ddef})$.
By \cref{lm:small-diverse}, we can choose $R_0\subset R$ such that 
\begin{itemize}
    \item $|R_0|\leq|Q_0|-1$,
    \item vertices in $R_0$ have different traces on $Q_0$ in $F$,
    \item $R_0$ is a witness of $A_0^{\ddef}$, that is vertices in $Q_0$ have different traces on $R_0$ in $F$.
\end{itemize}

 We are now ready to define sets $A',B',B''$ and number $n'$ as follows (see \cref{fig:few-defects-sets}):
    \[
         A'=A_0^{\pro}\cup A_0^{\ddef}\cup Q_0, \quad B'=B_0\cup R_0,\quad B''=B_0^{\pro}\cup B_0^{\ddef}\cup B_0^{\tree}\cup R_0, \quad n'=2(|\mathcal{P}_0|+|\mathcal{D}_0|).
    \]
In the rest of this section we verify that these sets and the number satisfy the three properties from \cref{thm:auto}.

\begin{figure}[htbp]
	\centering
	\begin{tikzpicture}
		% Set layout for the few-defects case of Theorem~\ref{thm:auto}.
% As in the earlier figures, named nodes are available for later annotations.

\node[
	containment/ambient,
	minimum width=.60\linewidth,
	minimum height=26mm
] (few-A) {};

% The same sample classes as in the many-defects figure.
\foreach \fraction in {0.12,0.20,0.34,0.66,0.80,0.88} {
	\draw[black!65, line width=.1pt]
		($(few-A.north west)!\fraction!(few-A.north east)$) --
		($(few-A.south west)!\fraction!(few-A.south east)$);
}

% A_0^def is away from the representative row.  A_0^pro and Q_0 lie in
% that row and may overlap; all three regions contribute to A'.
\node[
	containment/selected,
	draw=black!65,
	fill=FigRed,
	fill opacity=.18,
	minimum width=.312\linewidth,
	minimum height=10mm
] (few-Adef) at ($(few-A.center)+(0,2.5mm)$) {};

\node[
	containment/selected,
	draw=black!65,
	fill=FigRed,
	fill opacity=.18,
	minimum width=.26\linewidth,
	minimum height=7mm
] (few-Apro) at ($(few-A.center)+(-.085\linewidth,-7.5mm)$) {};

\node[
	containment/selected,
	draw=black!65,
	fill=FigRed,
	fill opacity=.18,
	minimum width=.312\linewidth,
	minimum height=7mm
] (few-Qzero) at ($(few-A.center)+(0,-7.5mm)$) {};

\node[
	containment/ambient,
	minimum width=.82\linewidth,
	minimum height=24mm,
	below=14mm of few-A
] (few-B) {};

\node[
	containment/ambient,
	minimum width=.19\linewidth,
	minimum height=23mm
] (few-R) at ($(few-B.center)+(-.215\linewidth,0)$) {};

\node[
	containment/ambient,
	minimum width=.43\linewidth,
	minimum height=23mm
] (few-Bgood) at ($(few-B.center)+(.105\linewidth,0)$) {};

\node[
	containment/core,
	draw=none,
	outer sep=0pt,
	fill=FigPoolFill,
	minimum width=.15\linewidth,
	minimum height=16mm
] (few-Rzero) at ($(few-R.center)+(0,-2.5mm)$) {};

% Equal blue and peach bands: R_0 belongs to both B' and B''.
% Draw the outline after the striped fill to keep the border continuous.
\path[containment/overlap, draw=black!65, line width=.1pt]
	(few-Rzero.north west) rectangle (few-Rzero.south east);

\node[
	containment/selected,
	draw=black!65,
	fill=FigPoolFill,
	minimum width=.37\linewidth,
	minimum height=16mm
] (few-Bzero) at ($(few-Bgood.center)+(0,-2.5mm)$) {};

% Back the union with white so single regions match the peach stripes.
% These three subsets of B_0 can overlap; their union is B'_0.
\node[
	containment/core,
	draw=none,
	fill=white,
	minimum width=.19\linewidth,
	minimum height=7mm
] (few-Bpro) at ($(few-Bzero.center)+(-6.5mm,2mm)$) {};

\node[
	containment/core,
	draw=none,
	fill=white,
	minimum width=.19\linewidth,
	minimum height=7mm
] (few-Bdef) at ($(few-Bzero.center)+(13.5mm,2mm)$) {};

\node[
	containment/core,
	draw=none,
	fill=white,
	minimum width=.19\linewidth,
	minimum height=8mm
] (few-Btree) at ($(few-Bzero.center)+(1.5mm,-3mm)$) {};

% Translucent peach fills make pairwise and triple overlaps darker.
\foreach \component in {few-Bpro,few-Bdef,few-Btree} {
	\path[fill=FigOrange, fill opacity=.24, rounded corners=2pt, draw=none]
		(\component.north west) rectangle (\component.south east);
}

% Draw all three contours after the fills so overlaps do not hide them.
\foreach \component in {few-Bpro,few-Bdef,few-Btree} {
	\draw[black!65, line width=.1pt, rounded corners=2pt]
		(\component.north west) rectangle (\component.south east);
}

\node[font=\small, anchor=east] at ($(few-A.west)+(-3mm,0)$) {$A$};
\node[font=\small, anchor=east] at ($(few-B.west)+(-3mm,0)$) {$B$};

\node[font=\fontsize{7.3pt}{8.8pt}\selectfont, anchor=north, inner sep=0pt] at
	($(few-A.north west)!0.06!(few-A.north east)+(0,-1mm)$) {$A_1$};
\node[font=\fontsize{7.3pt}{8.8pt}\selectfont, anchor=north, inner sep=0pt] at
	($(few-A.north west)!0.27!(few-A.north east)+(0,-1mm)$) {$A_i$};
\node[font=\fontsize{7.3pt}{8.8pt}\selectfont, anchor=north, inner sep=0pt] at
	($(few-A.north west)!0.73!(few-A.north east)+(0,-1mm)$) {$A_j$};
\node[font=\fontsize{7.3pt}{8.8pt}\selectfont, anchor=north, inner sep=0pt] at
	($(few-A.north west)!0.94!(few-A.north east)+(0,-1mm)$) {$A_k$};
\foreach \fraction in {0.16,0.50,0.84} {
	\node[font=\fontsize{7.3pt}{8.8pt}\selectfont, inner sep=0pt] at
		($(few-A.north west)!\fraction!(few-A.north east)+(0,-2.3mm)$)
		{$\cdots$};
}

\node[containment/label, font=\fontsize{7.3pt}{8.8pt}\selectfont] at
	($(few-Adef.north west)+(2mm,-2mm)$) {$A_0^{\ddef}$};
\node[containment/label, font=\fontsize{7.3pt}{8.8pt}\selectfont] at
	($(few-Apro.north west)+(1mm,-1mm)$) {$A_0^{\pro}$};
\node[font=\fontsize{7.3pt}{8.8pt}\selectfont, anchor=north east, inner sep=0pt] at
	($(few-Qzero.north east)+(-1mm,-1mm)$) {$Q_0$};

% Representatives in the four displayed classes; only the endpoint labels
% are shown to keep the overlapping A_0 regions uncluttered.
\node[circle, fill=black!70, inner sep=0pt, minimum size=1.2mm]
	(few-q-one) at
	($(few-A.south west)!0.06!(few-A.south east)+(0,3.5mm)$) {};
\node[font=\fontsize{7.3pt}{8.8pt}\selectfont, anchor=west, inner sep=0pt] at
	($(few-q-one.east)+(1mm,0)$) {$q_1$};
\node[circle, fill=black!70, inner sep=0pt, minimum size=1.2mm]
	(few-q-i) at
	($(few-A.south west)!0.27!(few-A.south east)+(0,3.5mm)$) {};
\node[circle, fill=black!70, inner sep=0pt, minimum size=1.2mm]
	(few-q-j) at
	($(few-A.south west)!0.73!(few-A.south east)+(0,3.5mm)$) {};
\node[circle, fill=black!70, inner sep=0pt, minimum size=1.2mm]
	(few-q-last) at
	($(few-A.south west)!0.94!(few-A.south east)+(0,3.5mm)$) {};
\node[font=\fontsize{7.3pt}{8.8pt}\selectfont, anchor=west, inner sep=0pt] at
	($(few-q-last.east)+(1mm,0)$) {$q_k$};

\node[font=\fontsize{7.3pt}{8.8pt}\selectfont, inner sep=0pt] at
	($(few-A.south)+(0,-5mm)$) {$A_0=A_0^{\pro}\cup A_0^{\ddef}$};

\node[containment/label, font=\fontsize{7.3pt}{8.8pt}\selectfont] at
	($(few-R.north west)+(2mm,-2mm)$) {$R$};
\node[containment/label, font=\fontsize{7.3pt}{8.8pt}\selectfont] at
	($(few-Rzero.north west)+(2mm,-2mm)$) {$R_0$};
\node[containment/label, font=\fontsize{7.3pt}{8.8pt}\selectfont] at
	($(few-Bgood.north west)+(2mm,-2mm)$) {$B_{\good}$};
\node[containment/label, font=\fontsize{7.3pt}{8.8pt}\selectfont] at
	($(few-Bzero.north west)+(2mm,-2mm)$) {$B_0$};
\node[font=\fontsize{7.3pt}{8.8pt}\selectfont, anchor=north west, inner sep=0pt] at
	($(few-Bpro.north west)+(2mm,-1.5mm)$) {$B_0^{\pro}$};
\node[font=\fontsize{7.3pt}{8.8pt}\selectfont, anchor=north east, inner sep=0pt] at
	($(few-Bdef.north east)+(-2mm,-1.5mm)$) {$B_0^{\ddef}$};
\node[font=\fontsize{7.3pt}{8.8pt}\selectfont, inner sep=0pt] at
	($(few-Btree.center)+(0,-1mm)$) {$B_0^{\tree}$};
	\end{tikzpicture}
	\caption{The choice of \setcaptionbox{FigRed!25}{A'=A_0^{\pro}\cup A_0^{\ddef}\cup Q_0}, \setcaptionbox{FigBlue!25}{B'=B_0\cup R_0}, and \setcaptionbox{FigOrange!25}{B''=B_0^{\pro}\cup B_0^{\ddef}\cup B_0^{\tree}\cup R_0} in the few-defects case.}
	\label{fig:few-defects-sets}
\end{figure}

\paragraph{Property (1).} Using the assumption that $C_0$ is sufficiently large, we obtain
$$
 |A'|\leq |A_0^{\pro}| + |A_0^{\ddef}| + |Q_0|\leq |A_0^{\pro}| + 2|A_0^{\ddef}| <2(|\mathcal{P}_0|+|\mathcal{D}_0|)=n',
$$
%and
$$
 |B'|=|B_0|+|R_0|\geq \frac{C^{1/3}(|\cP_0| + |\cD_0|)}{2} = \frac{C^{1/3}}{4}n' = \frac{C_0^{4/3}}{4}n' \geq C_0n',
$$
and 
$$
|B''|<|\mathcal{P}_0|+|\mathcal{D}_0|+|E(T)|+|A_0^{\ddef}|<3(|\mathcal{P}_0|+|\mathcal{D}_0|)=3n'/2 < 3n'.
$$
Recalling \cref{rem:B'-truncation}, these establish Property (1).

\paragraph{Property (2).}
As in \cref{sc:many}, since $B'\setminus B'' \subset B_0$, in order to show that $B'\setminus B''$ is $A'$-diverse in $F$, we show that $B_0$ is $A'$-diverse in $F$.
Indeed, let $u,v\in B_0$ be distinct. Then $P_u,P_v\in\mathcal{P}_0$ and $D_u,D_v\in\mathcal{D}_0$. By \cref{obs:profile-defect}, we have $(P_u,D_u)\neq (P_v,D_v)$. 
First, suppose that $P_u\neq P_v$. Then, since all sets $P\cap A_0^{\pro}$, $P\in\mathcal{P}_0$, are different, we have $P_u\cap A_0^{\pro}\neq P_v\cap A_0^{\pro}$, which implies $N_{A'}(u)\neq N_{A'}(v)$, as required.
Now, suppose that $P_u=P_v$, but $D_u\neq D_v$. Then, since all sets $D\cap A_0^{\ddef}$, $D\in\mathcal{D}_0$, are different, we have $D_u\cap A_0^{\ddef} \neq D_v\cap A_0^{\ddef}$, and therefore there exists $x$ in $(D_u \triangle D_v) \cap A_0^{\ddef}$. 
Since $P_u=P_v$, we have $\Sigma(u)=\Sigma(v)$, and thus $D_u\triangle D_v = (N_A(u) \triangle \Sigma(u)) \triangle (N_A(v) \triangle \Sigma(v)) = N_A(u) \triangle N_A(v)$.
Consequently, $x$ is adjacent to exactly one of $u$ and $v$, implying
$N_{A'}(u)\neq N_{A'}(v)$.

\paragraph{Property (3).}
Let $H:=\mathrm{Aut}(F[A',  B''])$ and $H_{\prodef}:=\mathrm{Aut}(T)$ act on $A'\cup B''$ and $\mathcal{P}_0\cup\mathcal{D}_0$, respectively. As in the previous section, we consider stabilisers
\[
    H_0:=H_{B^{\pro}_0,B^{\ddef}_0,B^{\tree}_0,A_0^{\pro},A_0^{\ddef},Q_0,R_0} \qquad \qquad  (H_{\prodef})_0:=(H_{\prodef})_{\mathcal{D}_0,\mathcal{P}_0}.
\]
We have that $|(H_{\prodef})_{0}|\leq|H_{\prodef}|\leq C^{n'/2}<C_0^{3n'}$. By \cref{lm:orbit-stabiliser}, we have 
$$
 |H| \leq 2^{7 (|A'| + |B''|)}\cdot |H_0| \leq 2^{18n'}\cdot |H_0|.
$$
Therefore, in order to establish property (3), it suffices to show that there exists an injective homomorphism $\varphi:\, H_0\hookrightarrow(H_{\prodef})_0$. Indeed, due to $C_0$ being large enough, it would imply the desired bound 
\[
    |H| \leq 2^{18n'} \cdot |H_0| < (2^{6} \cdot C_0)^{3n'}<C_0^{100n'}.
\]

Since automorphisms in $H_0$ preserve $B_0^{\pro}$ and $B_0^{\ddef}$, we can define the homomorphism $\varphi$ in the natural way as the projection on $\mathcal{P}_0\cup \mathcal{D}_0$: for $\sigma\in H_0$, we let $\varphi(\sigma)(P)=P_{\sigma(v(P))}$ for $P\in\mathcal{P}_0$ and $\varphi(\sigma)(D)=D_{\sigma(u(D))}$ for $D\in\mathcal{D}_0$. Then $\varphi$ is a homomorphism to a group of permutations acting on $\mathcal{P}_0\cup\mathcal{D}_0$. We then need to show that $\varphi(\sigma)\in (H_{\prodef})_0$ for every $\sigma\in H_0$, so that $\varphi$ is a homomorphism between groups of automorphisms, and that it has a trivial kernel and so it is injective. In order to establish triviality of the kernel, observe that, since $\varphi$ is the projection, any permutation that acts non-trivially on $B_0^{\pro}\cup B_0^{\ddef}$ maps to a non-trivial permutation. Therefore, it suffices to show that only the trivial permutation in $H_0$ fixes all vertices in $B_0^{\pro}\cup B_0^{\ddef}$. 
Hence, the following claim completes the proof of Property (3), and therefore the proof of~\cref{thm:auto}. 

\begin{claim}
$\varphi$ is a homomorphism and $(H_0)_{(B_0^{\pro}\cup B_0^{\ddef})}$ is trivial.
\end{claim}

\begin{proof}
We start from the following observation: different vertices in $B_0^{\pro}$ have different neighbourhoods in $A_0^{\pro}$. Indeed, these neighbourhoods are intersections of their profiles with $A_0^{\pro}$, while $A_0^{\pro}$ distinguishes all sets from $\mathcal{P}_0$ and each profile in $\mathcal{P}_0$ has a unique representative in $B_0^{\pro}$. Similarly, different vertices in $B_0^{\ddef}$ have different intersections of their defects with $A_0^{\ddef}$. 

Now, let $\sigma\in H_0$, and let $P\in\mathcal{P}_0$ and $D\in\mathcal{D}_0$ be adjacent in $T$. We need to prove that $\varphi(\sigma)(P)$ and $\varphi(\sigma)(D)$ are adjacent in $T$. Let $v=v(P)\in B_0^{\pro}$ and $u=u(D)\in B_0^{\ddef}$. Then there exists $t\in B_0^{\tree}$ such that $t$ corresponds to the edge $\{P,D\}$, that is $P_t=P=P_v$ and $D_t=D=D_u$. Also let $v'=v(P_{\sigma(t)})$ and $u'=u(D_{\sigma(t)})$. Since $\sigma(t)\in B^{\tree}_0$, we have that $\{P_{\sigma(t)},D_{\sigma(t)}\}\in E(T)$.

Since 
$$
 N(v)\cap A_0^{\pro}=P\cap A_0^{\pro}=N(t)\cap A_0^{\pro},
$$
we get that 
$$
P_{\sigma(v)}\cap A_0^{\pro}=N(\sigma(v))\cap A_0^{\pro}=N(\sigma(t))\cap A_0^{\pro}=P_{\sigma(t)}\cap A_0^{\pro}=P_{v'}\cap A_0^{\pro}.
$$
Since $v',\sigma(v)\in B_0^{\pro}$, we get $v'=\sigma(v)$.

Also, from the second part of \cref{lem:rep}, we derive
$$
 (D_u\cap A_0^{\ddef}=D_t\cap A_0^{\ddef})\Rightarrow (D_{\sigma(u)}\cap A_0^{\ddef}=D_{\sigma(t)}\cap A_0^{\ddef}).
$$
This implies
$$
D_{\sigma(u)}\cap A_0^{\ddef}=D_{\sigma(t)}\cap A_0^{\ddef}=D_{u'}\cap A_0^{\ddef}.
$$
Since $u',\sigma(u)\in B_0^{\ddef}$, we get  $u'=\sigma(u)$.

So, $\varphi(\sigma)(P)=P_{\sigma(v)}=P_{v'}=P_{\sigma(t)}$ and $\varphi(\sigma)(D)=D_{\sigma(u)}=D_{u'}=D_{\sigma(t)}$ are adjacent in $T$, as needed.

 It remains to show that $(H_0)_{(B_0^{\pro}\cup B_0^{\ddef})}$ is trivial. Assume $\sigma\in (H_0)_{(B_0^{\pro}\cup B_0^{\ddef})}$. Let $t\in B_0^{\tree}$. Let $v=v(P_t)$ and $u=u(D_t)$. We have just proved that $\sigma(v)=v(P_{\sigma(t)})$ and $\sigma(u)=u(D_{\sigma(t)})$. Therefore, $v(P_{\sigma(t)})=v=v(P_t)$ and $u(D_{\sigma(t)})=u=u(D_t)$, implying $t=\sigma(t)$. So, $\sigma$ fixes all vertices in $B_0^{\tree}$. Since all vertices $x\in A_0^{\pro}$ have different neighbourhoods $\{v:\,x\in P_v\}$ in $B_0^{\pro}$, we get that all these vertices are fixed as well. Similarly, all vertices of $A_0^{\ddef}$ are fixed since, due to the second part of \cref{lem:rep}, $\sigma$ preserves the relation $x\in D_u$, and since all vertices $x$ in $A_0^{\ddef}$ have different sets $\{u:\,x\in D_u\}$ in $B_0^{\ddef}$.  Due to the first part of \cref{lem:rep}, $\sigma$ also fixes all vertices in $Q_0$. Finally, since $R_0$ is $Q_0$-diverse, $\sigma$ fixes all vertices in $R_0$. Therefore, $\sigma$ is the identity. This completes the proof.
\end{proof}

\section{Neighbourhood Complexity of Slow Factorial Classes}
\label{sc:props-proofs}

In this section we prove our results on the neighbourhood complexity of classes which are at the bottom end of the factorial layer, just a little faster than small classes. 

Our first result, \cref{pro:lower}, follows from considering the hereditary closure of a sequence of sparse  random bipartite graphs with suitably chosen part sizes. The way in which the speed is controlled is similar to the disproof of the implicit graph conjecture \cite{HH22}.   
 
\lowerprop* 
\begin{proof}For $\varepsilon \in (0,1)$, let $\eta =\varepsilon/4$ and $\delta =\varepsilon^2/64$. Let $G_n$ be a bipartite graph on $V=(A,B)$ where $|A|=n$ and $|B|=\left\lfloor n^{1+\delta }\right\rfloor$  satisfying the following conditions:
	\begin{enumerate}
		\item\label{itm:2}  $|E(G_n[S])| \leq (1+\eta)|S| \text{ for all }S\subseteq V(G_n)\text{ with }|S|\leq n^{1-\eta}$;
	\item\label{itm:3} there exists an $A$-diverse set $W\subseteq B$ with  $|W| =\Omega(n^{1+\delta})$.    
	\end{enumerate}
   We now show such graphs exist. 
	\begin{claim}\label{clm:RG} Let $G_n$ be a random bipartite graph on $V=(A,B)$ where $|A|=n$ and $|B|=\lfloor n^{1+\delta }\rfloor $ and each edge between $A$ and $B$ is present independently with probability $1/n$. Then with high probability $G_n$ satisfies Conditions  \ref{itm:2} and \ref{itm:3}. 	
	\end{claim}
	\begin{pocd}{\ref{clm:RG}} 
		 For 	Condition \ref{itm:2}  as $|A|+|B| \leq 2n^{1+\delta}$, the probability there exists a set $S\subseteq V$ of size $s\leq n^{1-\eta}$ with at least $ \lceil(1+\eta)s\rceil$ edges is at most 
		\begin{align*}\binom{|A|+|B|}{s}  \binom{\binom{s}{2}}{\lceil(1+\eta)s\rceil}\frac{1}{n^{\lceil(1+\eta)s\rceil}}
        &\leq \left(\frac{2n^{1+\delta}e}{s} \right)^s\left(\frac{s^2e}{2\lceil(1+\eta)s\rceil n} \right)^{\lceil(1+\eta)s\rceil} \\
        &\leq \left( C n^{\delta-\eta }s^\eta   \right)^s \leq \left( C n^{\delta-\eta^2 }    \right)^s,\end{align*}
		for some constant $C>0$, where the final inequality uses $s\leq n^{1-\eta}$. Now,  as $\eta^2 = \varepsilon^2/16=4\delta>0$ and by the union bound,  Condition \ref{itm:2} fails with probability at most 
		\[\sum_{s=1}^{\lceil n^{1-\eta}\rceil}  \left( C n^{ \delta-\eta^2 }    \right)^s\leq \sum_{s=1}^{\lceil n^{1-\eta}\rceil}  \left( C n^{-3\delta}    \right)^s =o(1). \]

		 For Condition \ref{itm:3} the probability that a fixed $b\in B$ has exactly two neighbours is $\binom{n}{2}\frac{1}{n^2}(1-\frac{1}{n})^{n-2}$ which is bounded from below by $1/(2e)>0$. It follows from the Chernoff bound that with high probability there are $\Omega(|B|) = \Omega(n^{1+\delta})$ many vertices in $B$ with exactly two neighbours in $A$, call this set $B'$. The neighbours of a fixed $b\in B'$ are uniform over all pairs of vertices in $A$ and independent of other neighbourhoods of vertices in $B'$. Thus, conditional on $|B'|$, 
         the expected number of pairs $b_1,b_2 \in B'$ with the same neighbourhood is \[\binom{|B'|}{2}\cdot \binom{n}{2}^{-1} = O(n^{2\delta}) = o(|B|).\] Consequently, by Markov's inequality, with high probability there is a set $W\subset B'\subset B$ of size $\Omega(n^{1+\delta})$  with distinct neighbourhoods in $A$.  
	\end{pocd}  
	
	We now take the class $\mathcal{C} = \her(\bigcup_{n\geq 1} G_n)$ to be the hereditary closure of these graphs. Observe that, as $\delta= \varepsilon^2/64$,  Condition \ref{itm:3} establishes the lower bound on the neighbourhood complexity of $\mathcal{C}$. Thus, it remains to show that $|\mathcal{C}_n| =O((n!)^{1+\varepsilon})$.  

For an $n$-vertex graph $H\in \mathcal{C}$, let $N$ be the smallest integer such that $H$ is an induced subgraph of the graph $G_N$. We split the proof over two cases; $(i)$: $n \leq N^{1-\eta}$, and $(ii)$: $n >N^{1-\eta}$. 
\begin{labeling}{\textit{Case (ii)}:}
	\item[\textit{Case (i)}:] In this case there exists a set $S\subset V(G_N)$ of size  $n\leq N^{1-\eta}$ such that $H=G_N[S]$. Thus by Condition \ref{itm:2} we have $|E(H)| \leq (1+\eta)n  $. However, the number of $n$-vertex graphs with at most $(1+\eta)n  $ edges is at most  
	\[
 \sum_{j=0}^{\lfloor (1+\eta)n \rfloor}  \binom{\binom{n}{2}}{j } \leq 2n\cdot \binom{\binom{n}{2}}{\lfloor (1+\eta)n \rfloor } \leq 2n\cdot \left( \frac{n^2 e}{ (1+\eta)n} \right)^{(1+\eta)n}= 	2^{O(n)} \cdot (n!)^{1+\eta}.
	\] 
	\item[\textit{Case (ii)}:] We have that there are at most $n^{\frac{1}{1-\eta}}$ choices for the size $N$ for any graph $H$ in Case (ii). For each such $N$, the graph $G_N$ has at most $2N^{1+\delta}\leq 2n^{\frac{1+\delta}{1-\eta}}$ many vertices. Observe that each labelled induced subgraph $H$ from such a $G_N$ comes from an injective function $f:[n] \to V(G_N)$. As there are at most $|V(G_N)|^n$ such functions, the number of such graphs $H$ is bounded by \[n^{\frac{1}{1-\eta}} \cdot (2n^{\frac{1+\delta}{1-\eta}})^{n} =2^{O(n)} \cdot n^{\frac{1+\delta}{1-\eta}n} =2^{O(n)} \cdot n^{\frac{1+\varepsilon^2 /64}{1-\varepsilon/4}n} \leq  	2^{O(n)} \cdot (n!)^{1+\varepsilon/2}.   \]      
\end{labeling}
Thus combining the two bounds above gives
\[|\mathcal{C}_n| =  2^{O(n)} \cdot \left(  (n!)^{1+\eta} +  (n!)^{1+\varepsilon/2}\right)  = 2^{O(n)} \cdot(n!)^{1+\varepsilon/2} = O((n!)^{1+\varepsilon}),\]which establishes the bound on the speed of $\mathcal{C}$, and thus the result. 
\end{proof}

We now restate and prove \cref{pro:upper}, whose proof follows \cite[Theorem 4]{BDSZ25} closely. 
 \upperprop*
	
	\begin{proof}
	Let $\mathcal{C}$ be a~hereditary class of graphs such that $|\mathcal{C}_n| \leqslant (n!)^{1+\varepsilon} $ holds for every $n \in \mathbb{N}$.
	To prove the claim, we will show that for every sufficiently large $n$, every graph $G$ in $\mathcal{C}$ and every set $A\subseteq V(G)$ of size $n$ it holds that
	\[
	|\{N_G(v)\cap A~:~v \in V(G) \}|\leq  2n^{1+\delta} \]where $\delta= \varepsilon + (\log n)^{-1/2}$. 
	Suppose this is not the case. Then, by definition, there exists a~graph $G \in \mathcal{C}$ and an $n$-element set $A \subseteq V(G)$ such that 
	$|\{N_G(v)\cap A~:~v \in V(G) \}| > 2n^{1+\delta}$.
	Since vertices in $A$ contribute at most $n$ different neighbourhood traces on $A$, there exists a set $B\subseteq V(G)\setminus A$ such that $|B|>2n^{1+\delta}-n\ge n^{1+\delta} $ and $B$ is $A$-diverse in the bipartite graph consisting of the edges of $G$ between $A$ and $B$. Let $F=G[A,B]$.

	Let $n_1 := \lfloor n^{1+\gamma} \rfloor $ and $N := n_1 + n$, where $\gamma=\frac{(\log n)^{-1/2} }{2(1+\varepsilon)}$.
	We will show that the hereditary closure of $G$ contains more than $(N!)^{1+\varepsilon}$ labelled graphs on $N$ vertices, contradicting our assumption on the number of graphs in $\mathcal{C}$. 
	To do this, let us first fix a~labelling of $A$ with $[n]$ which will be the same for all our labelled graphs. Then, we construct the family of labelled graphs by taking subsets of $B$ of size $n_1$ and labelling each of them with the integers from $n+1$ to~$N$. Any two such induced graphs are distinct since $B$ is $A$-diverse, and restoring the edges within $A$ or $B$ cannot destroy this property. 
	
	Hence, the number of $N$-vertex labelled graphs in $\mathcal{C}$ that are induced subgraphs of $G$ is lower bounded by the number of ways to choose a~labelled subset of size $n_1$ from $B$, which is
	\begin{equation}\label{eq:lower-bound}
		n_1! \cdot \binom{|B|}{n_1} \geqslant 	n_1! \cdot \left(\frac{|B|}{n_1}\right)^{n_1} \geq  n_1! \cdot n^{(\delta -\gamma)n_1}.
	\end{equation}
	We want to show that this is more than the total number of $N$-vertex labelled graphs in $\mathcal{C}$, which, by assumption and as $\varepsilon, \gamma \leq 1$ and $n$ is large, does not exceed
	\begin{equation}\label{eq:upper-bound}
		(N!)^{1+\varepsilon}   = ((n_1+n)!)^{1+\varepsilon}  \leqslant (n_1!)^{1+\varepsilon} \cdot (n_1+n)^{(1+\varepsilon)n} \leq  n_1!  n^{(1+\gamma)\varepsilon n_1} \cdot n^{5n}.
	\end{equation}
	From \cref{eq:lower-bound} and \cref{eq:upper-bound}, we should then show that 
	\[
	n_1! \cdot n^{(\delta -\gamma)n_1} >  n_1!  n^{(1+\gamma)\varepsilon n_1} \cdot  n^{5n},
	\]
	which, by cancelling and taking logs, follows from establishing that
	\[\delta -\gamma > (1+\gamma) \varepsilon + \frac{5n}{\lfloor n^{1 +\gamma}\rfloor},	\]which is implied by 		\[\delta - \varepsilon - (1+\varepsilon)\gamma >  10n^{-\gamma} . 	\] 
	Observe that this holds for sufficiently large $n$ after we recall that $\varepsilon\in (0,1)$ is fixed, $\delta = \varepsilon + (\log n  )^{-1/2}$ and $\gamma=\frac{(\log n)^{-1/2} }{2(1+\varepsilon)}$. 
\end{proof}

\bigskip
\paragraph{Acknowledgments.}
This work was initiated during a Research-in-Groups stay at James Clerk Maxwell’s birthplace in Edinburgh. The authors are grateful to the International Centre for Mathematical Sciences for providing a stimulating and productive research environment. This work was supported by an EPSRC grant ``Hereditary graph classes: from Size to Structure''.  Viktor Zamaraev is supported by the Leverhulme Research Fellowship (RF-2026-309$\backslash$9).

\paragraph{Disclosure of AI assistance.} This project originated at ICMS in Edinburgh in April 2026. During that meeting, we reduced the main result to~\cref{thm:auto} and developed an initial proof strategy. This approach involved randomly selecting a pool $B_{\good}$ of vertices from $B$ and then using group-theoretic methods to inductively construct subsets $A'\subset A$ and $B'\subset B_{\good}$, as required by~\cref{thm:auto}. Although our principal group-theoretic result was correct, in late September 2026 we discovered a flaw in the induction step. We subsequently used ChatGPT Pro (Astra 6) to assist in repairing the argument. This led to a modification of the original approach in which the randomly selected pool was supplemented by an additional, deterministically constructed set of edges obtained via~\cref{lm:profile-defect}. Both the statement and a proof of this lemma were suggested by ChatGPT. In parallel, ChatGPT proposed an alternative strategy for constructing the sets $A'$, $B'$, and $B''$ that avoided both random selection and the group-theoretic machinery of our original approach. After independently verifying the validity of this strategy, we decided to adopt it in the present paper, as it yields a substantially shorter and less technical proof. Thus, the overall combinatorial strategy and the statements of the auxiliary results in~\cref{sec:proof-of-auto} were suggested by ChatGPT, while their proofs were developed independently by the authors. In particular, for~\cref{lm:profile-defect}, we developed our own, simpler probabilistic proof, which we present in this paper.

The proof of~\cref{pro:lower} was adapted from a similar construction suggested by ChatGPT Plus (Astra 6), which was based on subdividing the edges of a random graph.

ChatGPT was also used during the final stages of manuscript preparation to help identify typographical errors and inconsistencies. The entire manuscript has been written, carefully checked, and verified by the authors, who take full responsibility for the correctness of all statements, arguments, and results presented in the paper.

\bibliographystyle{alpha}
\bibliography{references}

\end{document}